\documentclass[11pt]{article}
\small
\normalsize
\usepackage[latin1]{inputenc}
\usepackage[T1]{fontenc}
\usepackage{fancyheadings}
\usepackage{amsmath}
\usepackage{amssymb}
\usepackage{amsfonts}
\usepackage{amsbsy}
\usepackage{amscd}
\usepackage{theorem}
\usepackage{epsf}
\usepackage{psfrag}
\usepackage{multicol}
\allowdisplaybreaks

\def\build#1_#2^#3{\mathrel{\mathop{\kern 0pt#1}\limits_{#2}^{#3}}}
\def\noi{{\noindent}}
\def\be{\begin{equation}}
\def\ee{\end{equation}}
\def\ba{\begin{eqnarray*}}
\def\ea{\end{eqnarray*}}

\def\E{{\bf E}}
\def\aa{{\bf a}}
\def\cqfd{ \hfill $\blacksquare$ }

\def\llb{[\hspace{-.10em} [ }
\def\rrb{ ] \hspace{-.10em}]}

\def\s{{\cal S}}

\def\der{{\rm d}}
\def\a{{\bf a}}
\def\g{{\cal G}}
\def\aal{{\bf \alpha}}

\def\f{{\cal F}}

\def\e{{\cal E}}
\def\t{{\cal T}}
\def\v{{\cal V}}

\def\p{{\cal P}}

\def\u{{\cal U}}
\def\P{{\bf P}}

\def\ee{\epsilon}
\def\xxe{x_{\epsilon}}
\def\zxe{\zeta_{x,\epsilon}}
\def\txe{\tau_{\epsilon}^{x}} 
\def\tbar{\overline{\tau}}
\def\Tbar{\overline{\cal T}}
\def\tbaree{\overline{\tau}_{\epsilon}}
\def\Tbaree{\overline{\cal T}_{\epsilon}}
\def\tbareex{\overline{\tau}_{\epsilon}^{x}}
\def\Tbareex{\overline{\cal T}_{\epsilon}^{x}}
\def\mubar{\overline{\mu}}
\def\trr{{\bf Tr}}
\def\ppp{{ p}}
\def\sh{{\bf Sh}}
\def\til{\widetilde{\tau}}
\def\Til{\widetilde{{\cal T}}}
\def\tilee{\widetilde{\tau}_{\epsilon}}
\def\Tilee{\widetilde{{\cal T}}_{\epsilon}}
\def\tileex{\widetilde{\tau}_{\epsilon}^{x}}
\def\tileexp{\widetilde{\tau}_{\epsilon_p}^{x}}
\def\Tileex{\widetilde{{\cal T}}_{\epsilon}^{x}}
\def\muti{\widetilde{\mu}}
\def\laweq{\overset{(law)}{=}}
\def\du{\frac{d}{u}}

\def\tilH{\widetilde{H}}
\def\tilC{\widetilde{C}}
\def\Ubar{\overline{U}}
\def\Vbar{\overline{V}}
\def\dist{{\bf d}}
\def\spin{{\bf Sp}}
\def\ztileex{\widetilde{\zeta}_{x, \ee}}

\def\l{{\cal L}}

\def\noi{\noindent}

\def\tte{\tau_{\epsilon}}

\def\build#1_#2^#3{\mathrel{\mathop{\kern 0pt#1}\limits_{#2}^{#3}}}

\newtheorem{theorem}{Theorem}[section]
\newtheorem{lemma}[theorem]{Lemma}
\newtheorem{proposition}[theorem]{Proposition}

{\theorembodyfont{\rmfamily}\newtheorem{remark}{Remark}[section]}
\newcommand{\R}{\mathbb{R}}
\newcommand{\RR}{\mathbb{R}}

\newcommand{\Z}{\mathbb{Z}}

\newcommand{\N}{\mathbb{N}}
\newcommand{\U}{\mathbb{U}}
\newcommand{\F}{\mathbb{F}}

\newcommand{\T}{\mathbb{T}}
\newcommand{\un}{\boldsymbol{1}}

\begin{document}

\title{ {\bf  CONTINUUM TREE LIMIT FOR THE RANGE \\
OF RANDOM WALKS ON REGULAR TREES} }
\author{ by \\
Thomas Duquesne,\\
{\small Université Paris 11, Mathématiques, 91405 Orsay Cedex, France. } \\
{\small e-mail: thomas.duquesne@math.u-psud.fr }}
\vspace{4mm}
\date{\today } 

\maketitle

\begin{abstract}
Let $b$ be an integer greater than $1$ and let $W^{\ee}=(W^{\ee}_n; n\geq 0)$
be a random walk on the $b$-ary rooted 
tree $\U_b$, starting at the root, 
going up (resp. down) with probability $1/2+\epsilon$ (resp. $1/2 -\epsilon$) , $\epsilon \in (0, 1/2)$, and choosing direction $i\in \{ 1, \ldots , b\}$ 
when going up with probability $a_i$. Here $\aa =(a_1, \ldots, a_b)$ 
stands for some non-degenerated fixed set of weights. We consider the range 
$\{ W^{\ee}_n ; n\geq 0 \}$ that is a subtree of $\U_b $. It corresponds to a
unique random rooted ordered tree that we denote by $\tau_{\epsilon}$. We rescale the edges of
$\tau_{\epsilon}$ by a factor $\ee $ and we let $\ee$ go to $0$: we prove 
that correlations due
to frequent backtracking of the random walk only give rise to a
deterministic phenomenon taken into account by a positive 
factor $\gamma (\aa)$. More
precisely, we prove that 
$\tau_{\epsilon}$ converges to a continuum random tree encoded by two independent Brownian motions with drift conditioned to
stay positive and scaled in time by $\gamma (\aa)$. We actually state 
the result in the more general case of a random walk on a tree 
with an infinite number of branches at each node ($b=\infty$) 
and for a general set of weights $\aa =(a_n, n\geq 0)$.  


\vspace{4mm}

\noi {\it MSC 2000 subject classifications:} 60F17, 60J80, 05C05, 05C80.

\vspace{4mm}

\noi {\it Key words and phrases:} continuum random tree, contour process,
exploration process, height process, limit theorem,  
random walk, range, regular tree.

\end{abstract}

\section{Introduction.}

  Random walks on trees have been intensively studied by many authors 
having different motivations coming from  
group theory, discrete potential theory, statistical mechanics or genetics. 
We refer to the book of W. Woess \cite{Woess} for a general introduction 
to random walks on infinite graphs and to the book of 
R. Lyons and Y. Peres \cite{LyPe} for a probabilistic approach more focused on trees. 
See also \cite{LyPemPer97} for a survey of open problems concerning 
random walks on trees. In most of the papers about random walks on trees,  
given the tree-like environment  
the transition probabilities of the random walk are fixed and one 
focuses on a certain range of questions: 
the speed of the random walk (see \cite{SaSt87} for random walks
on groups, \cite{Tak97} for random walks on periodic trees, \cite{LyPemPe95}
and \cite{LyPemPe96}, for random walks on Galton-Watson trees), large
deviation principle for the distance-from-the-root process 
(see \cite{GaPeZe} for random walks on Galton-Watson trees), 
central-limit theorem for the distance-from-the-root process and the number of
visited vertices (see \cite{ChYaXi} for the $b$-ary tree and \cite{Piau98} 
for the simple random walk on supercritical Galton-Watson trees). 
In this paper, we consider a different problem; the transition
probabilities are not fixed: we study, near 
criticality, transient random walks on the $b$-ary rooted tree and 
more generally on the $\infty$-ary tree, in 
a ``diffusive'' regime.

  Let us specify 
that we only consider ordered rooted trees that are formally defined
as in \cite{Ne}: Let $\N= \{ 0, 1, 2, \ldots \} $ be the set of 
the nonnegative integers, set 
$\N^* = \N \setminus \{ 0\}$. The $\infty$-ary tree is the set 
$\U = \{ \varnothing \} \cup \bigcup_{n\geq 1} (\N^* )^n $
of the finite words written with 
positive integers by. Let $u\in \U$ be the word $u_1 \ldots u_n$, $u_i \in
\N^{*}$. We denote the length of $u$ by $|u|$ : $|u|=n$. 
$|u|$ is viewed as {\it the height} of the vertex $u$ in $\U$. Let $v=v_1 \ldots v_m
\in \U $.  Then the word 
$uv$ stands for the concatenation of $u$ and $v$: $uv=u_1 \ldots u_n v_1
\ldots v_m $. Observe that $\U$ 
is totally ordered by the {\it  lexicographical
order} denoted by $\leq $ . A rooted ordered tree $t$ is a subset of $\U$ satisfying the following 
conditions
\begin{description}
\item{(i)} $\varnothing \in t$ and $\varnothing $ is called the {\it root} of $t$.

\item{(ii)} If $v\in t$ and if $v=uj$ for some $j\in \N^* $, then, $u\in t $. 

\item{(iii)} For every $u\in t $, there exists $k_u (t) \geq 0 $ such that $uj\in t$ for every $1\leq j\leq k_u (t)$. 
\end{description}
We denote by $\T$ the set of ordered rooted trees. Let
us mention that we sometimes see ordered rooted trees as family trees. So, 
we often use the genealogical terminology instead of 
the graph-theoretical one. 
All the random objects introduced in this paper are defined on 
an underlying probability space denoted by $(\Omega , \f , \P)$. Let $\ee \in
(0,1/2)$ and let $\aa=(a_n, n\geq 1)$
be some non-degenerated fixed set of weights, namely  
$\sum a_n =1$ and $0\leq a_n <1$ , $n\geq 1$.  

   We attach to the infinite tree $\U$ a cemetery point 
$\partial \notin \U$ situated at height $(-1)$ and we view $\partial$ as the 
parent of the root $\varnothing $. Then, we let run a particle on 
$\U \cup \{ \partial \}$ that evolves as follows:

\begin{itemize}

\item The particle starts at $\varnothing$ at time $0$ and it stops 
when it reaches  $\partial$. 

\item If at time $n$ the particle is at vertex $v\in \U$, 
then it jumps down to the parent of $v$ with probability 
$1/2-\epsilon$ and it goes up with probability $1/2+\epsilon$.

\item When going up, the particle chooses direction $j\in \N^*$ and jumps 
to the vertex $vj\in \U$ with probability $a_j$. 

\end{itemize}

   The height of the particle evolving in $\U \cup \{ \partial \}$ is then distributed as a random walk 
on $\Z$ started at $0$, stopped when reaching state $-1$, 
and whose possible jumps are $(+1)$ with probability 
$1/2+\epsilon$ and $(-1)$ with probability 
$1/2-\epsilon$.  
In this paper {\bf we condition the particle to 
never reach $\partial$} (observe that this conditioning is non singular). We denote by $W^{\ee}= (W_n^{\ee} \, ; n\geq 0)$ 
the sequence of vertices in $\U$ visited by the conditioned particle.

  We study the range $\{ W_n^{\ee} \, ; n\geq 0 \}$ when $\epsilon$ goes to
zero. Observe that it is an ordered rooted subtree of 
$\U $. There exists a unique ordered rooted tree $\tau_{\ee}
\in \T$ corresponding to $\{ W_n^{\ee }\, ; n\geq 0 \}$ via a one-to-one 
map that fixes the root $\varnothing$, preserves adjacency and that is increasing with respect to the lexicographical order. 

  Since $W^{\ee} $ goes to infinity, 
$\tau_{\ee}$ has one single infinite line of 
descent. Following Aldous's terminology introduced in \cite{Alfringe} we call
{\it sin-tree} such trees (see  Section 2.1 for precise definitions). The
distribution of $\tau_{\ee}$ is not simple and it 
shows correlations due to frequent backtraking of the random walk (see
comments in Section 2.3). However Theorem \ref{theolim}, 
which is the main result of the paper, asserts that $\tau_{\ee}$ converges 
in distribution to some continuum random tree. More precisely, think of $\tau_{\ee}$ as a planar
graph embedded in the clockwise oriented half-plane and suppose 
that its edges have length one; consider a particle visiting 
continuously the edges of 
$\tau_{\ee} $ at speed one from the left to the right, going backward as less as possible;
we denote by  $ C_s (\tau_{\ee} ) $ 
the distance from the root of the particle at time $s$
and we call the resulting process $C(\tau_{\ee})=(C_s (\tau_{\ee}) ; s\geq 0
)$ the {\it left contour 
process} of $\tau_{\ee}$. It is clear that the particle never reaches the part
of $\tau_{\ee}$ at the right hand of the infinite line of descent; observe however
that $C(\tau_{\ee})$ completely encodes the left part of $\tau_{\ee}$. Denote by 
$C^{\bullet}(\tau_{\ee})$ the process corresponding to a particle visiting 
$\tau_{\ee}$ from the right to the left. Thus, $(C(\tau_{\ee}), C^{\bullet}
(\tau_{\ee}))$ completely encodes 
$\tau_{\ee}$ (see Section 2.2 for more careful definitions and other encodings of sin-trees).  
Let $D$ and $D^{\bullet}$ be two independent copies of the process 
$s\rightarrow B_s -2s-2 \inf_{r\leq s} (B_r -2r)$ where $B$ is 
distributed as the standard linear Brownian motion started at $0$. 
Theorem \ref{theolim}
asserts that the following convergence 
$$\left( \ee C_{s/\ee^2}(\tau_{\ee}) \; , \;  \ee C^{\bullet}_{s/\ee^2} (\tau_{\ee}) 
\right)_{s\geq 0}
\xrightarrow[\ee\rightarrow 0]{\quad \quad} \left( 2D_{\gamma s} \;, \;2D^{\bullet}_{\gamma s}\right)_{s\geq 0}, $$
holds in distribution in $C([0,\infty), \R^2)$ endowed with the topology 
of uniform convergence on compact sets. 
We see that correlations in $\tau_{\ee}$ only give rise 
to a deterministic phenomenon characterized
by a constant $\gamma = \gamma (\aa)$ that is defined by 
\begin{equation}
\label{coeff}
1/\gamma = \E \left[ \left( 1+ X_1 + X_1 X_2 + X_1 X_2 X_3 + \ldots
  \right)^{-1} \right] , 
\end{equation}
where $(X_n \,; n\geq 1)$ stands for a sequence of i.i.d. 
$\{a_n, \, n\geq 1\}$-valued random variables whose distribution is given by 
$\P (X_n =a_i)= \sum a_j $, the sum being taken over the $j$'s such that
$a_j=a_i$. Observe that if $b $ is some integer greater than $1$ 
and if $a_n=0 $ for all
$n\geq b+1$, then 
the particle remains in the $b$-ary ordered rooted tree $\U_b = \{ \varnothing \} \cup \bigcup_{n\geq 1} \{ 1, \ldots , b \}^n $.
More comments about this limit theorem 
are added before and after the statement of Theorem \ref{theolim}. 

   Before ending this section, let us give a short overview of the proof of
the theorem: one part of the proof relies on a specific encoding of 
the range $\{ W_n^{\ee} \, ; n\geq 0 \}$ that can be 
explained as follows: Denote by $(|  W_n^{\ee}| \, ; n\geq 0)$ the sequence of
successive heights of the particle. It is obviously distributed as 
a random walk started at $0$ whose possible jumps are 
$(+1)$ with probability $1/2+\epsilon$ and $(-1)$ with probability 
$1/2-\epsilon$, conditioned to stay nonnegative. Then, the piecewise linear 
process
$$ t\longrightarrow |  W_{\lfloor t \rfloor }^{\ee}| + 
(t-\lfloor t \rfloor)|W_{\lfloor t \rfloor +1 }^{\ee}| $$
is the contour process of an infinite 
``fictive'' tree denoted by $\tbar_{\epsilon}$ whose distribution can be informally
described as follows: $\tbar_{\epsilon}$ has one
infinite line of descent; at each vertex $v$ on the infinite line of 
descent an independent random number with distribution $\mu$ 
of independent Galton-Watson trees with offspring distribution 
$\mu$ is attached at the left of the infinite line. Here, $\mu$ stands for the probability measure on $\N$ given by 
$\mu(k)=(1/2+\epsilon)(1/2-\epsilon)^k$, $k\geq 0$ (see 
Section 2 for precise definitions concerning trees and Lemma 3.1 for 
the details). 

  We then encode the walk $( W_n^{\ee} \, ; n\geq 0 )$ by the tree 
$\tbar_{\epsilon}$ and random marks $\mubar_u \in \N^*$ , 
$u\in \tbar_{\epsilon}$ that are defined as follows:  Let $u\in 
\tbar_{\epsilon}$ be 
distinct from the root $\varnothing $. Denote $\overleftarrow{u}$ its parent. 
By definition of the contour process 
the edge $(\overleftarrow{u}, u)$ corresponds to a unique upcrossing 
of the process $(|W_n^{\ee}| \, ; n\geq 0)$ between times $n(u)$ and 
$n(u)+1$. Thus, there exists $j\in \N^*$ such that the word 
$W_{n(u)+1}^{\ee}$ is written $W_{n(u)}^{\ee} \, j$ and we set 
$\mubar_u =j$. Then, we easily check that conditional on 
$\tbar_{\epsilon}$, the marks $\mubar_u $ , $u\in \tbar_{\epsilon}
\setminus{\{ \varnothing \}}$ are independent and 
distributed on $\N^*$ in accordance
with $\aa$ (see Section 3.1 for details). We get back the walk 
$W^{\ee}$ from the marked tree $\Tbar_{\epsilon} 
=(\tbar_{\epsilon}\; ;\;  (\mubar_u, u\in \tbar_{\epsilon}))$,  in the
following way: consider $u\in 
\tbar_{\epsilon}$, distinct from the root $\varnothing $ at height 
$|u|=n$; denote
by $u_0 = \varnothing $ , 
$ u_1$ , $\ldots$, $ u_n =u $ the ancestors of
$u$ listed in the genealogical order. Then we define the {\it track} of 
$u$: $ \trr_{\Tbar_{\epsilon}} (u)$ by the word 
$\mubar_{u_1} \ldots \mubar_{u_n} \in \U $ (observe that the mark of the root
plays no role). Then, 
$$ W_{n(u)+1}^{\ee}= \trr_{\Tbar_{\epsilon}} (u) $$
and thus 
$$ \trr_{\Tbar_{\epsilon}} (\tbar_{\epsilon})=\{ W_n^{\ee} \, ; n\geq 0 \} .$$
Taking the trace of $\tbar_{\epsilon}$ has two distinct effects: the first 
one shuffles $\tbar_{\epsilon}$ in the order of the marks in 
$\N^*$. The second one shrinks the tree because several edges of 
$\tbar_{\epsilon}$ might correspond to the same vertex in $\U$.

Let us briefly explains how to deal with the shuffling effect of the tree: it
is possible to reorder randomly the marked tree $\Tbar_{\epsilon}$ into a new 
marked tree $\Til_{\epsilon} = (\til_{\epsilon} \; ;\;  (\muti_u \, , u\in 
\til_{\epsilon}))$ such that:

\begin{itemize}

\item $\til_{\epsilon}$ has the same distribution as the tree obtained from 
$\tbar_{\epsilon}$ by changing independently and uniformly at random the 
order of birth of brothers in $\tbar_{\epsilon}$.

\item If $u_1$ , $u_2\in \til_{\epsilon}$ are 
such that $u_1\leq u_2$ then 
$$\trr_{\Til_{\epsilon}} (u_1) \leq \trr_{\Til_{\epsilon}} (u_2)  $$
\end{itemize}
(see Section 3.1 for a precise definition). 
Thus, the shuffled tree $\til_{\epsilon}$ has a simple 
distribution specified by Remark 3.1. Up to the shrinking effect, 
$\til_{\epsilon}$ is close to $\tau_{\epsilon}$ and if we denote by 
$(C(\til_{\epsilon}), C^{\bullet} (\til_{\epsilon}))$ the left and the right
contour processes of $\til_{\epsilon}$ we prove in Section 2.1 that 
\begin{equation}
\label{premconv}
\left( \ee C_{s/\ee^2}(\til_{\epsilon}) \; , \;  \ee 
C^{\bullet}_{s/\ee^2} (\til_{\epsilon}) 
\right)_{s\geq 0}
\xrightarrow[\ee\rightarrow 0]{\quad \quad} \left( 2D_s \;, \;2
D^{\bullet}_{s}\right)_{s\geq 0}, 
\end{equation}
in distribution in $C([0,\infty), \R^2)$ endowed with the topology 
of uniform convergence on compact sets. Denote by 
$\dist $ the graph distance in $\U$. Informally speaking, 
(\ref{premconv}) says that the metric space 
$(\til_{\epsilon} , \ee \, .\,  \dist)$ converges to some 
random metric space that is a continuum random trees encoded by  
$D$ and $D^{\bullet}$ (see Comment (c) after Theorem 2.1 for a more precise 
discussion of that point). Next, observe that for any 
$u_1$ and $u_2$ in $\til_{\epsilon}$ 
$$ 0\leq \dist (u_1,u_2) - 
\dist \left( \trr_{\Til_{\epsilon}} (u_1), \trr_{\Til_{\epsilon}} (u_2)  \right) \leq 2 G (u_1 , u_2)$$
where conditional on $u_1, u_2 \in \til_{\epsilon}$, 
$G (u_1, u_2)$ is a random integer with a geometric distribution with
parameter $q=\P(X_1 \neq X_2) $, which is a quantity that does not 
depend on $\ee$. 
This gives an informal argument to explain why the limits of 
the metric spaces $(\tau_{\epsilon} , \ee \, .\, \dist)$ and 
$(\til_{\epsilon} , \ee \, .\, \dist)$ should be close 
and why tightness for the contour processes of $\tau_{\epsilon}$ is not the 
difficult part of the proof: it is deduced from (\ref{premconv}) 
by (now standard) arguments 
inspired from the proof of Theorem 20 in Aldous's paper \cite{Al2}.

   The technical point of the paper concerns the identification of the limiting 
tree  by studying precisely the ``shrinking effect'' via explicit 
computations for the $\U$-indexed Markov process $(Z_v, v\in \U)$ given by 
$$ Z_v= \# \{ u \in \til_{\epsilon} \; : \; \trr_{\Til_{\epsilon}} (u)=v
\}.$$
This analysis is done in Proposition 3.3 and in Proposition 3.4.  
More precisely, if we fix a real number 
$x>0$ and if we remove from $ \til_{\epsilon}$ all 
the descendents of the unique vertex at height $\lfloor x/\ee \rfloor$ on the
infinite line of descent, we get a finite tree denoted by 
$\til_{\epsilon}^x$. Let $\cal{U}$ be a uniform random variable in 
$[0,1]$ independent of $\til_{\epsilon}$. Denote by $U(\ee)$ 
the vertex of $\til_{\epsilon}^x$ coming in 
the $\lfloor U \#\til_{\epsilon}^x \rfloor$-th position 
in the lexicographical order. Set 
$$ \Vbar (\ee)= \sum_{\substack{v\in \U  \\ v\leq 
\trr_{\Til_{\epsilon}} (U(\ee))}} \un_{\{ Z_v  >0 \}} . $$
Then, we prove that 
$$\ee^2 \left( \Vbar(\ee) - \frac{1}{\gamma} 
 {\cal U} \# \tilde{\tau}_{\epsilon}^x \right) 
\xrightarrow[\ee\rightarrow 0]{\quad  \quad} 0 $$ 
in probability. This key result is stated more precisely in Lemma 
\ref{convproba}.

   The paper is organized as follows: In Sections 2.1 and 2.2 we specify our 
notations 
and we define various encodings of trees and forests; 
Theorem \ref{theolim} is stated at Section 2.3; Section 3 is devoted to 
its proof that relies on a certain combinatorial representation of the range $\{ W_n^{\ee}; n\geq 0 \} $ given at Section 3.1 
and on a technical estimate (Lemma \ref{convproba}) whose proof is postponed at
Section 3.3 while the proof of 
Theorem \ref{theolim} itself is done at Section 3.2.

\section{Preliminaries and definitions}

\subsection{Trees, forests and sin-trees}

We first start with some notations. We define on $\U$ the {\it genealogical}
order $\preccurlyeq $ by
$$\forall u,v \in \U \; , \quad u\preccurlyeq v \Longleftrightarrow \exists w\in \U \; :\; v=uw  .$$
If $u\preccurlyeq v$, we say that $u$ is an ancestor of $v$. If $u$ is distinct from the root,
it has an unique predecessor with respect to $\preccurlyeq$ that 
is called its parent and that is denoted by $\overleftarrow{u}$. 
We define the youngest common ancestor of $u$ and $v$ by the 
$\preccurlyeq$-maximal element $w\in \U$ such that $w\preccurlyeq u $ and $ w \preccurlyeq v $ and we denote it by 
$u\wedge v$. We also define the distance between $u$ and $v$ by $\dist (u,v)=|u| +|v|-2 |u\wedge v|$ and we use notation 
$\llb u,v \rrb$ for the shortest path 
between $u$ and $v$. Let $t\in \T$ and $u\in t$. We define the tree $t$
shifted at $u$ by 
$\theta_u (t)= \{ v\in \U \; :\; uv\in t\}$ and  
we denote by $[t]_u $ the tree $t$ cut at the node $u$  :  $[t]_u := \{ u\} \cup 
\{ v\in t \; :\; v\wedge u \neq u \}$. Observe that $[t]_u \in \T$. For any
$u_1, \ldots ,u_k \in t$ we also set 
$[t]_{u_1, \ldots , u_k} : =[t]_{u_1} \cap \ldots \cap [t]_{u_k}$ and
$$ [t]_n= [t]_{\{u\in t:\; |u|=n   \}}=\{ u\in t \; :\; |u| \leq n \} \; ,
\quad n\geq 0. $$

Let us denote by $\g $ the $\sigma $-field on $\T$ generated by the sets 
$\{ t\in \T :\;  u\in t\}$ , $u\in \U$ and let 
$\mu $ be a probability distribution on $\N$. We call {\it Galton-Watson tree} with offspring distribution 
$\mu$ (a GW($\mu $)-tree for short) any 
$(\f , \g)$-measurable random variable $\tau $ 
whose distribution is characterized by the two following 
conditions: 
\begin{description}
\item{(i)} $\P (k_{\varnothing }(\tau )=i) = \mu (i) \; ,\; i\geq 0 $.

\item{(ii)} For every $i\geq 1$ such that $\mu (i) \neq 0$, the shifted trees $\theta_{1} (\tau) , \ldots , \theta_{i} (\tau) $  
under \\ 
$\P (\cdot \mid k_{\varnothing } (\tau )=i)$ are independent copies of 
$\tau$ under $\P$.
\end{description}

\begin{remark}
\label{indGW}
Let $u_1, \ldots , u_k \in \U$ such that $u_i \wedge u_j \notin \{ u_1, \ldots
, u_k\}$ , $1\leq i,j\leq k$, and 
let $\tau $ be a GW($\mu $)-tree. Then, conditional on the event 
$\{ u_1, \ldots , u_k \in \tau\} $, $\theta_{u_1} (\tau) , \ldots , \theta_{u_k} (\tau) $ are i.i.d. GW($\mu $)-trees independent
of $[\tau]_{u_1, \ldots , u_k}$. 
\end{remark}

  We often consider a forest (i.e. a sequence of trees) instead of a single tree. More precisely, we define 
the forest $f$ associated with the sequence of trees $(t_l\, ; l\geq 1)$ by
the set 
$$ f= \{ (-1,\varnothing)\} \cup \bigcup_{l\geq 1} \left\{ (l,u) ,\; u\in t_l \right\} $$
and we denote by $\F $ the set of forests. Vertex $ (-1,\varnothing)
$ is viewed as a fictive root situated at generation $-1$. Let $u'=(l,u) \in f$ with $l\geq
1$; the height of $u'$ is defined by $|u'| :=|u|$ and its ancestor is defined
by $(l, \varnothing )$. For convenience, we denote it by 
$\varnothing_l:=(l, \varnothing )$. As already
specified, all the ancestors $\varnothing_1, \varnothing_2$, ...  are the 
descendants of $  (-1,\varnothing) $ and are situated at generation $0$. 
Most of the notations concerning trees extend to forests: 
The lexicographical order $\leq $ 
is defined on $f$ by taking first the individuals of $t_1$, next those of
$t_2$ ... etc and leaving $  (-1,\varnothing) $ unordered. The genealogical order $\preccurlyeq$ 
on $f$ is defined tree by tree in an obvious way. Let $v'\in f$. The youngest common ancestor of $u'$ and $v'$ 
is then defined as the $\preccurlyeq$-maximal element of $w'$ such that  $w'\preccurlyeq u'$ and $w'\preccurlyeq v'$ and 
we keep denoting it by $u'\wedge v'$. The number of children of $u'$ is $k_{u'} (f):=k_u (t_l)$ and the 
forest $f$ shifted at $u'$ is defined as the tree $\theta_{u'} (f) :=\theta_{u} (t_l)$. We also define 
$[f]_{u'}$ as the forest $\{ u'\} \cup \{ v'\in f \; :\; v'\wedge u' \neq u' \}$ and we extend 
notations $[f]_{u_1', \dots , u_k'} $  and $[f]_n$ in an obvious way. For
convenience of notation, we often identify
$f$ with the sequence $(t_l \, ; l\geq 1)$. When   $(t_l \, ; l\geq 1)=(t_1,
\ldots , t_k, \varnothing , \varnothing , \ldots )$, we say that $f$ is a
finite forest with $k$ elements and we abusively write $f=(t_1, \ldots ,
t_k)$.

We define the set of sin-tree by 
$$ \T_{sin}= \left\{  t\in \T \; :\; \forall n \geq 0 , \quad  \# \left\{ v\in t \; :\; |v|=n \; {\rm and } \;  \# \theta_v (t) =\infty 
\right\} =1 \; \right\} .$$
Let $t\in \T_{sin}$. For any $n\geq 0$, we denote by $u^*_n (t)$ the unique
individual $u$ on the infinite line of descent (i.e. such that $\# \theta_{u}
(t) =\infty $) situated at height $n$. Observe that
$u^*_0 (t)=\varnothing$. We use notation $\ell_{\infty} (t)=\{ u^*_n (t) ; \; n\geq
0 \}$ for the infinite 
line of descent of $t$ and we denote by $(l_n (t) \, ; n\geq 1)$ the sequence of positive integers such that
$u^*_n (t)$ is the word $l_1 (t) \ldots l_n (t) \in \U$.
We also introduce the set of {\it sin-forests} $\F_{sin}$ that is defined as the set of forests 
$f=(t_l \, ; l\geq 1)$ such that all the trees $t_l$ are finite except one
sin-tree $t_{l_0}$. We extend to sin-forests notations 
$u^*_n$, and $l_n$ by setting $l_n (f)=l_n (t_{l_0})$, $u^*_n (f)=(\, l_0, u^*_n (t_{l_0})\, )$ and $u^*_0 (f)=\varnothing_{l_0}$.

  Next, we introduce a natural class of random sin-trees called {\it Galton-Watson trees with immigration} 
(GWI-trees for short). 
The distribution of a GWI-tree is characterized by

\noindent
$\bullet$ its {\it offspring distribution} $\mu $
on $\N$ that we suppose critical or subcritical: $\bar{\mu}=\sum_{k\geq 0} 
k\mu (k) \leq 1$;

\noindent
$\bullet$ its {\it dispatching distribution} $r$ defined on the first octant 
$\{ (k,l)\in \N^* \times \N^* \; :\; 1\leq l \leq k \}$ that prescribes the
distribution of the number of 
immigrants and their positions with respect to the infinite
line of descent.

\noindent
More precisely, $\tau$ is a GWI($\mu,r$)-tree if it satisfies the two following conditions:
\begin{description}
\item{(i)} The sequence $S=( \; (k_{u^*_n (\tau)} (\tau ), l_{n+1} (\tau)) \; ; \; n\geq 0)$ is i.i.d. with distribution $r$.

\item{(ii)} Conditional on $S$, the trees $\theta_{u^*_n (\tau)i} \, (\tau)$
  with $n\in \N $ and $1\leq i\leq k_{u^*_n (\tau)} (\tau )$ with $i\neq
  l_{n+1} (\tau)$ are mutually independent  GW($\mu$)-trees.
\end{description}

We define a GWI($\mu,r$)-forest with $l\geq 1$ 
elements by the  forest $\varphi=(\tau , \tau_1 , \ldots ,
\tau_{l-1}) $ 
where the $\tau_i$'s are i.i.d. GW($\mu$)-trees 
independent of the  GWI($\mu,r$)-tree $\tau$. It will be sometimes 
convenient to insert $\tau$ at random in the sequence 
$(\tau_1 , \ldots , \tau_{l-1})$ but unless otherwise specified
the random sin-tree in a random sin-forest occupies the
first row.

  The word ``immigration'' comes from the following obvious observation: 
Let $\varphi$ be a GWI($\mu,r$)-forest with $l+1$
elements. Set for any $n \geq 0$ , $Z_n (\varphi )=\# \{ u\in \varphi \; : \;
|u|=n \}-1$. Then the process $( Z_n (\varphi ); n\geq 0)$ 
is a {\it Galton-Watson process with immigration} started at state $l$, 
with offspring distribution $\mu $ and {\it immigration distribution} $\nu $
given by 
$$\nu (k)= \sum_{1 \leq j\leq k+1} r(k+1,j) \; , \; \quad k\geq 0. $$
Recall that 
a Galton-Watson process with immigration $( Z_n (\varphi ); n\geq 0)$ 
is a $\N$-valued Markov chain whose transition probabilities 
are characterized by 
\begin{equation}
\label{transGWI}
\E \left[ x^{Z_{n+m} (\varphi )} \mid Z_m (\varphi )\right]
= f_n (x)^{Z_{m} (\varphi )} g\left( f_{n-1} (x) \right)g\left( f_{n-2} (x) \right) \ldots g(f_0(x)),
\end{equation} 
where $f$ (resp. $g$) stands for the generating function of $\mu $
(resp. $\nu$) and 
where $f_n$ is recursively defined  by $f_n= f_{n-1} \circ f$ , $n\geq 1$ and
$f_0 ={\rm Id}$.

We conclude this section by giving an elementary result on the so called 
GW($\mu$)-{\it size-biased trees} that are GWI($\mu,r$)-trees
with dispatching distribution of the form $r(k,j)=\mu (k)/\bar{\mu} $, $1\leq j\leq k$. Size-biased trees
arise naturally by conditioning critical or subcritical GW-trees 
on non-extinction: see \cite{AlPit98}, \cite{Alfringe}, 
\cite{Kes87} or \cite{LyoPemPer} for related results. 
The term ``size-biased'' can be justified by the following 
elementary result needed at Section 3.3: Let $\varphi$ be a random
forest corresponding to a sequence of $l$ independent GW($\mu$)-trees and let  
$\varphi_{\flat}$ be a GWI($\mu,r$)-forest with $l$ elements where $r$ is taken as above and  
where the position of the unique random sin-tree in $\varphi_{\flat}$ is 
picked uniformly at random among the $l$
  possible choices. Check 
that for any nonnegative
measurable functional $G$ on $\F \times \U$:
\begin{equation}
\label{sizebias}
\E \left[ 
\sum_{u\in \varphi } G\left( [\varphi]_u ,  u \right)\right] = \sum_{n\geq 0} l \, \bar{\mu }^n  \E \left[ 
 G\left( [\varphi_{\flat}]_{ u^*_n (\varphi_{\flat})} , u^*_n (\varphi_{\flat})\right) \right]
\end{equation} 
and in particular ${\rm d}\P ( [\varphi_{\flat}]_n \in \cdot \, )/{\rm d}\P ( [\varphi]_n \in \cdot \, ) =  Z_n (\varphi )/ l\bar{\mu }^n$.

\subsection{The encoding of sin-trees}

The purpose of the paper is to provide a limit theorem for $\tau_{\ee}$ 
thanks to its encoding by two contour processes as 
briefly explained at the introduction. It will be convenient to
introduce two additional encoding processes: namely the {\it height process}
(also called {\it exploration process}) and a certain kind of random walk.

\noindent {\bf Encoding of finite trees and forests.} 
Let $t\in \T$ be a finite tree and let $u_0= \varnothing < u_1 < \ldots <
u_{\# t-1} $ be 
the vertices of $t$ listed in the lexicographical order. We define the {\it height process} of $t$ by $ H_n (t) = |u_n| $, 
$0\leq n<\# t$. $H(t)$ clearly characterizes the tree $t$.

  We also encode $t$ by its contour process which is 
informally defined as follows: think of $t$ as a graph embedded in the
clockwise oriented half-plane with unit length edges; 
let run a particle starting at the root at time $0$ that explores $t$ from the left to the right moving continuously along each edge 
at unit speed until it comes back to its starting point. In this evolution, each edge is crossed twice and the total amount of 
time needed to explore the tree is thus $2(\# t-1)$. The contour process 
$C(t)=(C_s(t); \, 0\leq s\leq 2(\# t -1) )$ is 
defined as the distance-from-the-root process 
of the particle at time $s\in [0\, ,\,2(\# t-1) ]$. 
More precisely, $C(t)$ can be recovered from 
the height process by the following transform: Set $b_n =2n-H_n (t)$ for $0\leq n <\# t$ and $b_{\# t}= 2(\# t-1)$. Then 
observe that
\begin{equation}
\label{contourvsheight}
C_s (t)= \left\{   
\begin{array}{lll}
\displaystyle 
H_n (t)-s+b_n &{\rm if} \; \displaystyle  s\in [b_n , b_{n+1}-1) \; {\rm and } \; n<\# t-1  \\
\displaystyle s-b_{n+1} +H_{n+1} (t) &{\rm if} \displaystyle \;  s\in [b_{n+1}-1 ,b_{n+1} ]\; {\rm and } \; n<\# t-1 \;  \\
\displaystyle H_{\# t-1} (t) - s+b_{\# t -1} &{\rm if} \displaystyle  \; s\in [b_{\# t-1} ,b_{\# t} ] .
\end{array} 
\right. 
\end{equation}

  We also need to encode $t$ in a third way by a path $V(t)=(V_n (t); 0\leq
n\leq \# t )$ that is defined by $  V_{n+1} (t) =
V_n (t) + k_{u_n} (t) -1$ and $V_0 (t)=0$. $V(t)$ is sometimes called the 
Lukaciewicz path associated with $t$. 
It is clear that we can reconstruct $t$ 
from $V(t)$. Observe that the jumps of $V(t)$ are $\geq -1$. Moreover 
$V_n (t)\geq 0$  for any $0\leq n < \#t $ and $V_{\# t} (t)=-1$. 
We recall from \cite{LGLJ1} without proof the following formula that allows
to write the height process as a functional of $V(t)$:
\begin{equation}
\label{heightvswalk}
H_n (t)= \# \left\{ 0\leq j<n \; :\; V_j (t)= \inf_{j\leq k\leq n } V_k (t) \right\} \;  ,\; 0\leq n <\#t. 
\end{equation}

\begin{remark}
\label{harrisrem}
If $\tau $ is a critical or subcritical GW($\mu $)-tree, then it is clear from our definition that $V(\tau)$ 
is a random walk started at $0$ that is stopped at $-1$ and whose jump
distribution is given by $\rho (k)= \mu (k+1)$ , $ k\geq -1$.  
However neither $H(\tau)$ nor $C(\tau)$ are Markov processes except for the
geometric case: $\mu(k)=(1-p)p^k$ with $p\in (0, 1/2]$. In this case, $C(\varphi)$ is
distributed as a random walk killed at $-1$ and  
whose possible jumps are $(+1)$ with probability $p$ and $(-1)$ 
with probability $1-p$ (more precisely 
it is the restriction of the first $T_{-1}-1$ 
steps of a random walks killed at
the reaching time of level $(-1)$).
\end{remark}

The previous definition of $V$ and of the height process can be easily
extended to a forest $f=(t_l \, ; l\geq 1)$ of finite trees 
as follows: Since all the trees $t_l$ are finite, it is possible to list 
all the vertices of $f$ but $(-1,\varnothing)$ in the lexicographical order:
$u_0 =\varnothing_1 < u_1 < \ldots $ etc. We then simply 
define the height process of $f$ by $H_n (f)=|u_n|$ and $V(f)$ by 
$V_{n+1} (f)=V_n(f) +k_{u_n} (f) -1$ with $V_0 (f)=0$. Set $n_p = \# t_1 +
\ldots +\# t_p $ and $n_0=0$ and observe that 
$$ H_{n_p+k} (f)= H_k
(t_{p+1}) \quad {\rm and} \quad
V_{n_p+k} (f)= V_k (t_{p+1})-p \; , \quad  
0\leq k<\# t_{p+1} \; , \, p\geq 0.$$ 
We thus see that the height 
process of $f$ is the concatenation of the height processes of the trees 
composing $f$. Moreover the $n$-th visited vertex $u_n$ is in $t_p$ iff
$p=1-\inf_{0\leq k\leq n} V_k (f)$. Then, it is 
easy to check that (\ref{heightvswalk}) remains true for every $n\geq 0$ when $H(t)$ and $V(t)$ are replaced by resp. 
$H(f)$ and $V(f)$.


\vspace{5mm}

\noindent {\bf Encodings of sin-trees.} Let $t\in \T_{sin}$. A particle visiting $t$ in the lexicographical order 
never reaches the part of $t$ at the right hand 
of the infinite line of descent. So we need two height processes or equivalently 
two contour processes to encode $t$. 
More precisely, the left part of $t$ is the set $\{ u\in t:\; \exists v\in
\ell_{\infty } (t)\; {\rm s.t.} \; u\leq v\}$. It can be listed in a
lexicographically increasing sequence of individuals denoted by $\varnothing = u_0<u_1< \ldots $.
We simply define the {\it left height process} of $t$ 
by $H_n(t)=|u_n|$ , $n\geq 0$. $H(t)$ completely encodes the left part of $t$. 
To encode the right part we consider the ``mirror
image'' $t^{\bullet}$ of $t$. More precisely, let $v\in t$ be the word
$c_1c_2 \ldots c_n $. For any 
$j\leq n$, denote by $v_j:=c_1 \ldots c_j$ the $j$-th ancestor of $v$ with 
$v_0=\varnothing$. Set 
$c^{\bullet}_j=k_{v_{j-1}} (t) -c_j +1$ and $v^{\bullet}=c^{\bullet}_1 \ldots
c^{\bullet}_n$. We then define $t^{\bullet}$ as 
$\{ v^{\bullet} , \;  v\in t \}$ and we define the {\it right height process} of $t$ as 
$ H^{\bullet}(t) : =H(t^{\bullet})$.
\begin{remark}
\label{symmetry}
Observe that $\tau$ and $\tau^{\bullet}$ have the same distribution if $\tau$
is a GW($\mu$)-tree. This is 
not anymore the case if $\tau$ is a  GWI($\mu, r$)-tree unless $r(k,m) = r(k, k-m+1)$.
\end{remark}
We now give a decomposition of $H(t)$ and $H^{\bullet}(t)$ along
$\ell_{\infty}(t)$ that is well suited to GWI-trees and that is
used at Section 3.2: Recall that $(u_n; n\geq 0)$ stands for the sequence 
of vertices of the left part of $t$ listed in the lexicographical 
order. Let us consider the set $\{ u_{n-1}^* (t)i ; 1 \leq i< l_n (t) ; n\geq 1\}$
of individuals at the left hand of $\ell_{\infty}(t)$ having a brother
on $\ell_{\infty}(t)$. To avoid trivialities, 
we assume that this set is not empty and we denote by $v_1< v_2 < \ldots $
the (possibly finite) sequence of its elements listed in the 
lexicographical order.

\begin{figure}[ht]
\psfrag{racine}{$0$}
\psfrag{premier}{$1\; ,\; 0*$}
\psfrag{deux}{$2\; ,\; 1*$}
\psfrag{trois}{$3\; ,\; 2*$}
\psfrag{quatre}{$4\; ,\; 3*$}
\psfrag{cinq}{$5\; ,\; 4*$}
\psfrag{six}{$6$}
\psfrag{sept}{$7$}
\psfrag{huit}{$8\; ,\; 5*$}
\psfrag{neuf}{$9\; ,\; 6*$}
\psfrag{dix}{$10\; ,\; 7*$}
\psfrag{onze}{$11\; ,\; 8*$}
\psfrag{douze}{$12$}
\psfrag{treize}{$13\; ,\; 9*$}
\psfrag{quatorze}{$14\; ,\; 10*$}
\psfrag{quinze}{$15$}
\epsfxsize=12cm
\centerline{\epsfbox{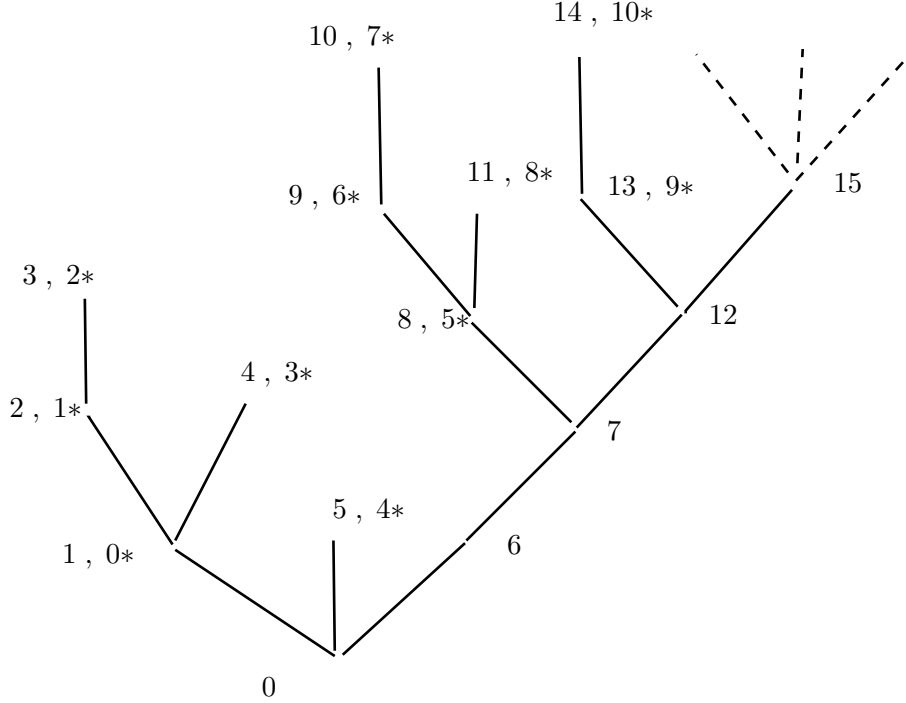}}
\caption{{\small The left part of a sin-tree $t$. The individuals 
which are not on $\ell_{\infty} (t)$ have two labels: 
the first one is their row in the lexicographical
order on $t$ and the second one (taged with a star) 
corresponds to their row in $f(t)$; individuals
of $\ell_{\infty} (t)$ have only one label corresponding to their row in $t$.}}
\label{leftsin.fig}
\end{figure}

  The forest $f(t)=(\theta_{v_1} (t) , \theta_{v_2} (t) , \ldots)$ is then
composed of the bushes rooted at the left hand of $\ell_{\infty} (t)$ taken in
the lexicographical order of their roots.
Define $L_n (t):=(l_1 (t)-1) + \ldots + (l_n(t)-1)$ , $n\geq 1$ with $L_0
(t)=0$ and consider the $p$-th individual of 
$f(t)$ with respect to the lexicographical order on $f(t)$; check that the
corresponding bush is rooted in $t$ at height 
$$ \aal (p)=\inf \{ k\geq 0 :\; L_k(t) \geq 1- \inf_{j\leq p } V_j (\ f(t)\ )\} \; .$$ 
Thus the corresponding 
individual in $t$ is $u_{{\bf n} (p)}$ where ${\bf n} (p)$ is given by 
\begin{equation}
\label{spinaldec1}
{\bf n}(p)=p+\aal (p)
\end{equation}
(note that the first individual of $f(t)$ is labelled by $0$). Conversely, 
let us consider $u_n$ that is the $n$-th individual of the left part of
$t$ with respect to the lexicographical order on $t$. Set 
${\bf p}(n) =\# \{ k<n : \; u_k\notin \ell_{\infty}(t)\}$ that is the number
of individuals coming before $u_n$ 
and not belonging to $\ell_{\infty}(t)$. Then
\begin{equation}
\label{spinaldec2}
{\bf p}(n)= \inf \{ p\geq 0\; :\; {\bf n}(p)\geq n \}. 
\end{equation}
and the desired decomposition follows:
\begin{equation}
\label{spinaldec3}
H_n (t)=n- {\bf p}(n)  + H_{{\bf p}(n)} (\, f(t) \, ).
\end{equation}
Since $n - {\bf p} (n) = \# \{ 0\leq k <n \; : \; u_k\in \ell_{\infty } (t)
\}$, we also get 
\begin{equation}
\label{spinaldec4}
\aal ( {\bf p} (n) -1) \leq n- {\bf p} (n) \leq \aal ({\bf p} (n)).
\end{equation}
Observe that if $u_n \notin \ell_{\infty} (t)$, then 
$n- {\bf p} (n)=\aal ({\bf p} (n))$. The proofs of these identities follow 
from simple counting arguments and they are left to the reader (see 
Figure (\ref{leftsin.fig})). 
Similar formulas
hold for $H^{\bullet}(t)$ taking 
$t^{\bullet}$ instead of $t$ in (\ref{spinaldec1}), (\ref{spinaldec2}),
(\ref{spinaldec3}) and (\ref{spinaldec4}).

\begin{remark}
\label{spinalGWI} 
The latter decomposition is particularly usefull when we consider a GWI($\mu ,r$)-tree $\tau$: In this case 
$(f(\tau), f(\tau^{\bullet}))$ is independent of $(L(\tau),
L(\tau^{\bullet}))$, $f(\tau)$ and $f(\tau^{\bullet})$
are mutually independent and 
$f(\tau)$ (resp.  $f(\tau^{\bullet})$) is a forest of i.i.d. GW($\mu
$)-trees if for some $k\geq 2$ we have $r(k,2)+\ldots 
+r(k, k)\neq 0$ (resp. $r(k,k-1)+\ldots 
+r(k, 1)\neq 0$), it is otherwise an empty forest. Moreover, the 
process $(L(\tau), L(\tau^{\bullet}))$ is a $\N \times \N$-valued random walk whose jump distribution is given by
$$ \P \left( L_{n+1} (\tau) -L_n (\tau)= m \, ; \,  L_{n+1} (\tau^{\bullet}) -L_n (\tau^{\bullet})=m' \right)= r(m+m'+1, m+1).$$
\end{remark}

We next define the left contour process of the sin-tree $t$ denoted
by $C(t)$ as the distance-from-the-root process of a particle starting at the root and moving  
clockwise on $t$ viewed as a planar graph embedded in the oriented half plane with edges of unit length. We define 
$C^{\bullet}(t)$ as the contour 
process corresponding to the anticlockwise journey and we can also write $C(t^{\bullet})=C^{\bullet}(t)$. More precisely, 
$C(t)$ (resp. $C^{\bullet}(t)$) can be recovered from  $H(t)$ (resp. $H^{\bullet}(t)$) through (\ref{contourvsheight})
that still holds for sin-trees (note that in that case 
the sequence $(b_n; n\geq 0)$ is infinite).

It will be sometimes convenient to approximate a sin-tree $t$ by the finite
tree $[t]_{u_n^*(t)}$ with $n$ large. The formula connecting the contour processes of $t$ and 
$[t]_{u_n^*(t)}$ is given as follows: 
Set $\sigma_n (t) = \# \{ u \in t : \; u < u_n^*(t) \} $
and $\sigma_n (t^{\bullet})=\# \{ u \in t^{\bullet} : \; u < u_n^*(t^{\bullet}) \}$.
We get 
$\sigma_n (t)+\sigma_n (t^{\bullet})= \# [t]_{u_n^*(t)} + n-1$ since the individuals of 
$\llb \varnothing , u_{n-1}^*(t)\rrb $ have been counted twice. Check that 
\begin{equation}
\label{cutcontour1}
\sigma_n (t)=\sup \{ k\geq 0 :\; H_k(t) \leq n \} \; \; {\rm and } \; \;  
2\sigma_n (t) -n =\sup \{ s\geq 0 :\; C_s(t) \leq n \},  
\end{equation}
with similar formulas for $t^{\bullet}$. Thus we get
\begin{equation}
\label{cutcontour2}
\left\{   
\begin{array}{ll}
\displaystyle 
C_s (t)= C_s ([t]_{u_n^*(t)})
&{\rm if} \; \displaystyle 
s\in [0,2\sigma_n (t) -n ] \; ,  \\
\displaystyle 
C^{\bullet}_s (t) = C_{ 2(\# [t]_{u_n^*(t)}-1) -s} \, ([t]_{u_n^*(t)})
&{\rm if} 
\displaystyle \; 
s\in [0,2\sigma_n (t^{\bullet}) -n ].
\end{array} 
\right. 
\end{equation}
(Observe that a similar formula is not available for height processes.)

\subsection{Statement of the main result}

For convenience of notation, we set $d=1/2 -\ee$ and $u=1/2 +\ee$. 
Recall that $\tte \in \T $ denotes the random ordered rooted tree associated with the range 
of the random walk $W^{\ee}$ in $\U$. First observe that the 
process $(|W^{\ee}_n|; n\geq 0)$ giving the distance from 
the root of the particle performing the random walk 
does contain an important part 
of the information concerning $\tte$. Moreover, this process is simply
distributed as the post-infimum path of a 
random walk whose possible jumps are $+1$ with probability $u$ and $-1$ with probability $d$. 
Recall that 
$(B_s\, ; s\geq 0)$ stands for the linear Brownian motion and set
for any $y\in \R$, 
$B_s^{(y)}=B_s+y s$ and $I_s^{(y)}=\inf_{u\leq s}B_u^{(y)}$.
Standard arguments imply 
$$ \left(\ee \left| W^{\ee}_{\lfloor s/\ee^2 \rfloor} \right|; s\geq 0\right) 
\xrightarrow[\ee\rightarrow 0]{\; ({\rm d})\;} \left( B_{s+ g }^{(2)}-I_{\infty}^{(2)} ; s\geq 0 \right), $$
where we have set $g= \inf \{ s\geq 0 : \;
B_{s}^{(2)}=I_{\infty}^{(2)}\}$. Notation 
$\overset{({\rm d})}{\rightarrow}$ stands for the convergence in distribution
in the appropriate space of right-continuous functions with left limits
endowed with Skorohod topology. We also use notation 
$\overset{({\rm fd })}{\rightarrow}$ 
for the convergence in distribution of all finite dimensional marginals.

  This result turns out to provide the right scaling for $\tte$ though the connection between 
$(|W^{\ee}_n| \, ; n\geq 0)$ and $\tte$ is non-trivial and the distribution of
$\tte$ is not simple: for instance, 
we can check that $\tte$ and $\tte^{\bullet}$ might not have 
the same distribution. Take
the binary case ${\bf a}=(a,1-a, 0, 0 , \ldots )$ for some $a\in (0,1)$. 
Define the set $A\subset \T$ by 
$A=\{ t\in \T : \; k_{\varnothing} (t)=2 , \; k_1 (t)=0, 
\; k_2(t) >0 \}$. Then it follows from simple arguments discussed in Section
3.1 that 
$$ \P (\tte \in A)= \frac{du^2a(1-a)}{(u+da)(u+d^2a)} \quad {\rm and } \quad 
\P (\tte^{\bullet} \in A)= \frac{du^2a(1-a)}{(u+d(1-a))(u+d^2(1-a))}.  $$
Thus, except for $a=1/2$, $\P (\tte \in A) \neq \P (\tte^{\bullet} \in A)$. 
Actually, when $\ee$ goes to zero, the 
particle backtracks more and more often causing correlations. However, Theorem \ref{theolim} asserts that
the correlations only give 
rise to a deterministic phenomenon that is taken into account by
the coefficient $ \gamma =\gamma (\aa)$ given by (\ref{coeff}).
\begin{theorem}
\label{theolim} 
Let $D$ and $ D^{\bullet}$ be two independent copies of $B^{(-2)}-2I^{(-2)}$. Then,
$$ (i) \quad \quad \left( \ee C_{s/\ee^2}(\tau_{\ee}) \; , \; \ee C^{\bullet}_{s/\ee^2} (\tau_{\ee}) 
\right)_{s\geq 0}
\xrightarrow[\ee\rightarrow 0]{\quad ({\rm d}) \quad} \left( 2D_{\gamma s}\; ,\;  2D^{\bullet}_{\gamma s}\right)_{s\geq 0} $$
$$ (ii)  \quad \quad \left( \ee H_{\lfloor s/2\ee^2 \rfloor }(\tau_{\ee}) \; ,
  \; \ee H^{\bullet}_{\lfloor 
s/2\ee^2 \rfloor} (\tau_{\ee}) 
\right)_{s\geq 0}
\xrightarrow[\ee\rightarrow 0]{\quad ({\rm d}) \quad} \left( 2D_{\gamma s} \; , \; 2D^{\bullet}_{\gamma s}\right)_{s\geq 0}. $$
\end{theorem}

\noindent Let us make some comments. {\bf (a)} The limit of the height and the
contour processes are the same up to the multiplicative
time constant $2$. This comes from the fact that vertices are visited
once by the height process while the edges are crossed exactly twice 
by the
contour process.

{\bf (b)} The definition of $\gamma $ through
 expectation (\ref{coeff}) is only for pratical reasons.
We have not found a simpler expression except for the case  
$a_1= \ldots = a_b = 1/b $ where $b$ is an integer greater than $1$. In that 
case the $X_i$'s are deterministic and 
$\gamma =1-1/b$.

{\bf (c)} The continuum random sin-tree whose $ 2D_{(\gamma \, \cdot )}$ and 
$ 2D^{\bullet}_{(\gamma \, \cdot)}$ are resp. 
the left and the right height processes can be defined as follows: To any real $s$ corresponds a vertex in the tree at height 
$H_s= \un_{(-\infty , 0)} (s) \, 2 D_{-\gamma s} +  \un_{[0, \infty )}
 (s) \, 2D^{\bullet}_{\gamma s}$ . Let $s\leq s'$. The youngest common ancestor of 
the vertices corresponding to $s$ and $s'$ is situated at height 
$$m(s, s')=
 \inf\{H_u ; \, u\in I(s,s') \} , $$ 
where $I(s, s')$ is taken as $[s,s']$ if $0\notin [s, s']$ and as 
$\R \setminus [s, s'] $ otherwise. Thus, the distance between the vertices correponding to $s$ and $s'$ is
$$\dist (s, s')= H_s + H_{s'}- 2 m(s, s') .$$
We say that $s$ and $s'$ are equivalent if they correspond to the same vertex in the tree, i.e.
$\dist (s,s')=0$ that is denoted by $s \thicksim s'$. We formally define the
 continuum random sin-tree as the quotient set 
$ T= \R / \thicksim $. Then $\dist $ induces a metric on $T$ that makes it be
a (random) Polish space.

We can show that the metric space $ (T ,\dist )$ is a $\RR$-tree
(see \cite{DuLG2} for related results). Due to the
Brownian nature of $H$, all fractal dimensions of $T$ are a.s. equal to $2$. 
A point $\sigma \in T$ is said a branching point if the open set $T\setminus
\{\sigma \} $ has more than two connected components and it corresponds to 
times at which $H$ reaches a local minimum. Since all the local minima of
$H$ are distinct, all the branching points are binary, i.e. $T\setminus
\{\sigma \} $ has $three$ connected components.


{\bf (d)} Observe that the limiting tree $T$ is symmetric since $D$ and
$D^{\bullet}$ have the same distribution. An heuristic explanation is the
following: arguments discussed at Section 3.2 imply that 
an unbalanced set of weights $\a $ breaks the symmetry of $\tte $
only if $\tte $ has branching points of order $\geq 3$ 
which does not happen to the limiting tree $T$ that is binary. 

\section{Proof of the main result}

\subsection{Combinatorial results}

In this section $\ee $ is fixed and for convenience of notation
we drop the corresponding subscript in the random
variables. Thus, we write $W$ and $\tau $ instead of $W^{\ee}$ and $\tte $. As
explained in the introduction, the linear interpolation of the process $(|W_n |\, ; n\geq 0)$ can be viewed 
as the left contour process of a (fictive) GWI-tree denoted by
$\tbar$ and whose distribution is given by the following lemma.
\begin{lemma}
\label{firsttree}
The linear interpolation of $(|W_n |\, ; n\geq 0)$ is distributed as the left contour process
of a GWI($\mu , r$)-tree 
where $\mu (k)=ud^k$, $r(k, k)=\mu (k-1)$ and $r(k,m)=0$, $1\leq m <k$ , $k\geq 0$.
\end{lemma}
\noindent 
{\bf Proof :} Let $(\xi_n; n\geq 1)$ be i.i.d. 
such that $\P(\xi_n=1)=u$ and $\P (\xi_n =-1)=d$. Set $S_0=0$ and 
$S_n=\xi_1+ \ldots + \xi_n $ and define $T_{-1}$ as $T_{-1}:= \inf \{ n\geq 0 \; :\; S_n=-1
\}$ (with the convention $\inf \emptyset =\infty$). Since the random walk 
$S=(S_n; n\geq 0)$ a.s. drifts to $+\infty$ , $P( T_{-1}=\infty ) >0$. By definition of $W$, 
$(|W_n |\, ; n\geq 0)$ has the same distribution as $S$ under 
$\P (\; \cdot \; |T_{-1}=\infty )$.  

Let us denote by $(T^{(0)}_i ; i\geq 0)$ the passage times to state $ 0 $: $T^{(0)}_0 =0$ and $ T^{(0)}_{i+1}= \inf\{n>T^{(0)}_i \; : \; S_n=0
\} $, with the convention $\inf \emptyset =\infty$. Set 
$$ K= \sup \left\{ i\geq 0 \; : \;  T^{(0)}_i < \infty \right\} 
< \infty  \quad {\rm
  a.s.} $$
We denote by $\e_1, \ldots ,\e_K, \e_{K+1}$ the excursions of $S$ away from $ 0$ defined by 
$$ \e_i=\left(S_{ T^{(0)}_{i-1} +n} ; 0\leq n\leq \zeta_i:=  T^{(0)}_{i}-
  T^{(0)}_{i-1} \right) \; ,\; 1\leq i\leq K  $$
and by $\e_{K+1}=(S_{ T^{(0)}_{K} +n} ; n\geq 0) $. We first consider the tree
$\tau$ whose contour process is the linear interpolation of $(S_n ; \; 0\leq n
\leq T_{-1} -1)$ under $\P (\; \cdot \; |T_{-1}<\infty )$:

\vspace{3mm}

{\bf Claim}:  $\tau$ is a GW($\mu$)-tree. 
\vspace{3mm}

\noindent 
{\bf Proof:} If $L$ stands for the number of children of the
root of $\tau$, then $L+1$ 
is also the number of times $S$ visits $ 0$ before $T_{-1}$: $L=\sup \{ i\geq 0 \; : \; T^{(0)}_{i} <T_{-1}
\}$. By applying the Markov property at the stopping times $ T^{(0)}_{i}$'s,
we show that $\P (L=0)=d$ and that for any $l\geq 1$, conditional on the event $\{ L=l\; ; \; T_{-1}<\infty\}$:  
\begin{itemize}
\item{(a)} $\e_1, \ldots , \e_{l}$ are i.i.d. and they are distributed as
  $\e_1$ under $\P ( \; \cdot \; | \e_1 (1) =1 \, ; \, T^{(0)}_{1} < \infty)$. 
\end{itemize}
\noindent
Moreover, the Markov property at time $1$ implies that 
\begin{itemize}
\item{(b)} $ ( \e_1 (n+1)-1 \; ; \; 0\leq n \leq T^{(0)}_{1} -2)$ under $\P (
 \;  \cdot \;  | \e_1 (1) =1 \, ; \, T^{(0)}_{1} < \infty)$ has the same distribution as $(S_n ;
  0\leq n \leq T_{-1} -1)$ under $\P (\; \cdot \; |T_{-1}<\infty )$. 
\end{itemize}
\noindent
Now observe that the contour processes of the subtrees $\theta_1 \tau , \ldots
, \theta_L \tau $ are the linear interpolations of 
$ ( \e_i (n+1)-1 \; ; \; 0\leq n \leq \zeta_i  -2)$ , $1\leq i\leq L$. We
deduce from (a) and (b) that $\tau $ satisfies the two conditions of the
definition of a GW-tree; its distribution is then the distribution 
of $L$ under 
$\P (\; \cdot \; |T_{-1}<\infty )$, which can be 
computed as follows: Observe first
that
$$ \{ L=l\; ; \; T_{-1}<\infty\}= \{  \e_1 (1)=1; T^{(0)}_{1} < \infty 
\; ; \; \ldots \; ; \;   \e_l (1)=1; T^{(0)}_{l} < \infty \; ; \; \e_{l+1}(1)=-1 \} .$$
Then,  by (a):
$$\P ( L=l \, ; \,  T_{-1}<\infty)= d\, 
\P (\e_1 (1)=1\, ;\,  T^{(0)}_{1} <\infty)^k . $$ 
But (b) implies that  $\P (\e_1 (1)=1; T^{(0)}_{1}<\infty )=u \, \P (  T_{-1}<\infty)$. Thus, 
$$ \P (L=l\; ; \; T_{-1}<\infty)= d  \left( u \, \P (  T_{-1}<\infty) \right) ^k  $$
and by summing over $l$ we get $\P (  T_{-1}<\infty)= d/ (1-u \P (
T_{-1}<\infty))$ which implies that $\P (  T_{-1}<\infty)=d/u$. Finally we 
get  
$$ \P ( L=l | T_{-1}<\infty)= ud^k =\mu (k) , $$
which achieves the proof of the claim. \cqfd 

\vspace{3mm}

Let us achieve the proof of the lemma: we now consider the tree $\tbar$ whose
contour process is the linear interpolation of 
$(|W_n |\, ; n\geq 0)$. To simplify notations, we identify this process to $S$ under 
$\P (\; \cdot \; |T_{-1}=\infty )$. Then $K+1$ is the
number of children of the ancestor of $\tbar$. First, observe that 
\begin{equation}
\label{setevent}
 \{ K=k \; ; \; T_{-1}=\infty\} = \{  \e_1 (1) =1; T^{(0)}_{1} < \infty 
\; ; \; \ldots \; ; \;   \e_k (1) =1; T^{(0)}_{k} < \infty \; ; \;  T^{(0)}_{k+1}
= \infty \} .
\end{equation}
By applying the Markov property, we
then show that conditional on $\{ K=k \; ; \; T_{-1}=\infty\}$ 
\begin{itemize}
\item $\e_1, \ldots , \e_{k+1}$ are independent;  

\item $\e_1, \ldots , \e_{k}$ are distributed as $\e_1$ under $\P ( \; \cdot
  \; |
  \e_1 (1)=1; T^{(0)}_{1} < \infty)$; 

\item $\e_{k+1}$ is  distributed as $\e_1$ under $\P ( \; \cdot \; | T^{(0)}_{1}
  =\infty)$
\end{itemize}
Now, by applying the Markov property at time $1$ we see that 
\begin{equation}
\label{tinfid}
 ( \e_1 (n+1)-1  ; \; n \geq 0 ) \; {\rm under} \; 
\P (\; \cdot \; | T^{(0)}_{1} = \infty)  \overset{(law)}{=} S  \; 
{\rm under} \; \P( \; \cdot \; |T_{-1}=\infty ) .
\end{equation}
Observe that the contour processes of the subtrees $\theta_1 \tbar , \ldots
, \theta_{K+1} \tbar $ are the linear interpolations of 
$ ( \e_i (n+1)-1 \; ; \; 0\leq n \leq \zeta_i  -2)$ , $1\leq i\leq
K+1$. Deduce from (b), from the previous claim and from (\ref{tinfid}) that 
conditional on $\{ K=k \; ; \; T_{-1}=\infty\}$,
the subtrees $\theta_1 \tbar , \ldots , \theta_{k} \tbar $ are $k$ independent
GW($\mu$)-trees and that $\theta_{k+1} \tbar $ is distributed as
$\tbar$. It implies that $\tbar$ satisfies the two conditions of the
definition of $GWI$-trees. Since the infinite subtree is $\theta_{k+1} \tbar $, $\tbar $ is a GWI($\mu , r$)-tree with 
$$r(k+1,m)=0 \; , \,  1\leq m <k+1 \quad 
 {\rm and } \quad 
r(k+1, k+1)=\P ( K=k |  T_{-1}=\infty ) \; , \; k\geq 0 , $$
which can be computed as follows: Deduce from (\ref{setevent}) and the Markov property 
$$ \P ( K=k \; ; \;   T_{-1}=\infty )= \P (\e_1 (1)=1\ ; \; T^{(0)}_{1} < \infty)^k \P ( T^{(0)}_{1} = \infty) . $$
Now observe that $\P ( T^{(0)}_{1} = \infty)=u\P( T_{-1}=\infty )$ and that 
$\P (\e_1 (1)=1\ ; \; T^{(0)}_{1} < \infty)= d$. Thus 
$\P ( K=k |  T_{-1}=\infty )=\mu (k)$, $ k\geq 0$, which achieves the proof of the lemma. \cqfd

\vspace{3mm}

Observe that $\tbar$ is completely asymmetric i.e. it has no vertices at the
right hand of its infinite line of descent. Note also that its immigration distribution $\nu $ is equal 
to $\mu $. In what follows, we explain how to recover the full
range $\{ W_n \, ; \, n\geq 0\}$ from $\tbar$. To that end 
we need to label $\tbar$ by random marks in $\N^*$ as explained in the
introduction. Let us introduce some notation: the set $T=(t \,; (m_u, u\in t))$ is a $\N^* $-marked tree $T$
if $t\in \T $ and if $m_u \in \N^* $, $u\in t$. The $m_u$'s are the 
marks of $T$. The set of $\N^* $-marked trees is denoted by
$\T_{\N^*} $. We define the
{\it track} of $T$ as the mapping  
$\trr_T\, :\, t \rightarrow \U $ defined as follows: Let $u\in t$; if we denote
by $u_0 = \varnothing 
\preccurlyeq u_1 \preccurlyeq \ldots \preccurlyeq u_n =u $ the ancestors of
$u$, then we define $\trr_T (u)$ as the word 
$m_{u_1} \ldots m_{u_n} \in \U $, with the convention $\trr_T (\varnothing )=
\varnothing$ (observe 
that $m_{\varnothing }$ plays no role in the definition of $\trr_T$).

  Similarly we define marked forets as sets of the form $F= (f; (m_u , u \in f ))$ where $f \in \F$ and 
$m_u \in \N^*$. The set of marked forests is denoted by 
$\F_{\N^*}$. We define the track $\trr_F$ of $F$ exactly as we have defined the
track of marked trees and we set for any $u\in F$
$$ \theta_u (F ) = \left(  \theta_u (f ) \, ; \,  (m_{uv} , v\in \theta_u (f)
  ) \right) 
\quad {\rm and } \quad [\, F \,]_u =\left( [\, f \,]_u \, ; \,
  (m_v , v\in [\, f \,]_u) \right).$$

Since the linear interpolation of the process $(|W_n |\, ; n\geq 0)$ is the
distance-from-the-root process of a (fictive) particle exploring continuously $\tbar$
at unit speed from left to right $\tbar$, we can associate 
with each vertex $u\in \tbar
\setminus\{ \varnothing \}$ a unique time $n(u) \in \N$ such that
the (fictive) particle climbs the edge 
$(\overleftarrow{u}, u)$ between times $n(u)$ and
$n(u)+1$. Since 
$|W_{n(u)+1}|=1+|W_{n(u)}| $, we can find $\mubar_u\in \N^*$
such that the word $W_{n(u)+1}$ is written $W_{n(u)}\mubar_u\in  \U$. 
We then define the
random marked tree $\Tbar$ as 
$$\Tbar =(\tbar\; ;\;  (\mubar_u, u\in \tbar)), $$
where the mark of the root $\mubar_{\varnothing }$ is taken independent of $W$
and distributed on $\N^* $ 
in accordance with the set of weights $\a$: 
$ \P (\mubar_{\varnothing }=i ) =a_i$, $i\in \N^*$
 The distribution of $\Tbar$  is described
by the elementary lemma whose proof is left to the reader. 
\begin{lemma}
\label{taubarlaw}
Conditional on $\tbar$, the marks $(\mubar_u, u\in \tbar)$ are independent and
distributed in accordance with $\a$. Moreover, 
\begin{equation}
\label{traceW}
\trr_{\Tbar} (\tbar) =\{ W_n \; ; \;  n\geq 0 \} .
\end{equation}
\end{lemma}

As already explained in the introduction, to take the track of $\tbar $ is a procedure that can be broken up in two
distinct sub-procedures:
The first one ``shuffles'' $\tbar$ by putting its edges 
in a certain random order. The second one ``shrinks'' $\tbar$ by identifying
some successive edges with respect to the new random
order. Let us first specify what we mean by {\it shuffling}: 
Let $t\in \T$; we say that $\ppp=( \ppp_u , u\in t)$ is a permutation of $t$ if
each $\ppp_u $ is a permutation of the (possibly empty) 
set $\{ 1, \ldots , k_u (t)\}$. Let $u\in t$ be the 
word $c_1 \ldots c_n$. We denote by $u_k=c_1 \ldots c_k$
the $k$-th ancestor of $u$. 
We define the word $u^{\ppp}$ by $\ppp_{u_0}(c_1) \ldots p_{u_{n-1}}(c_n)
\in \U$ if 
$u\neq \varnothing$ and by $\varnothing $ otherwise. We set 
$t^{\ppp}= \{ u^{\ppp} ; \, u\in t \}$. Now, pick uniformly at random a permutation 
$\pi $ of $t$ among the $\prod_{u\in t } k_u (t) ! $ possible ones. We define {\it the shuffling of } $t$ as the random tree $\sh (t):=t^{\pi}$.

\begin{remark}
\label{shuffleGWI} Shuffling a GW-tree does not change its distribution. It is also easy to check that $\sh (\tbar)$ is a GWI($\mu , r'$)-tree with $r'$ given by 
$r' (k,j)= ud^{k-1} /k$, $1\leq j \leq k$. 
\end{remark}

We would like to shuffle a $\N^*$-marked tree $T=(t \, ; (m_u ,\in t ))$ in accordance with the order of its marks 
in $\N^*$: for any permutation $p$ of $t$, set 
$T^{\ppp}= (t^{\ppp} \, ; (m_{u^{\ppp}} , u\in t))$ and observe that
\begin{equation}
\label{tra}
\trr_{T^{\ppp}}(t^{\ppp})= \trr_{T}(t).
\end{equation}
 Let $\pi(T)=(\pi_u , u\in t)$ 
be a random permutation of $t$ such that the $\pi_u$'s are mutually
independent and $\pi_u$ is picked uniformly at random among the permutations $\sigma$ of
$\{1, \ldots , k_u (t) \}$ satisfying 
$$m_{u\sigma(1)} \leq m_{u \sigma (2)} \leq \ldots \leq m_{u \sigma (k_u (t))} .$$ 
We define the shuffling of $T$ as $\sh (T):=T^{\pi (T)}$. By
definition the mapping 
$\trr_{\sh(T)}:t^{\pi (T)} \rightarrow \U $ is increasing with respect to the lexicographical order:
\begin{equation}
\label{increasing}
\forall u,v\in t^{\pi(T)} \; , \quad u\leq v \; \Longrightarrow \; \trr_{\sh(T)} (u) \leq \trr_{\sh(T)} (v) .
\end{equation}
Observe that if any brothers in $T$ have distinct marks, then $\pi (T)$ is
deterministic. Thus, $t^{\pi(T)}$ has clearly not the same distribution as 
$\sh (t)$. However when the
marks $m_u$, $u\in t$ are i.i.d. random variables, we can easily check that $t^{\pi (T)}$ is
distributed as $\sh (t)$. Thus, if we set 
$$\Til= \sh(\Tbar): = (\til \; ;\;  (\muti_u \, , u\in \til)) ,$$ 
then, we deduce from the previous observation that  
\begin{equation}
\label{ttilde1}
\trr_{\Til}(\til)=\{ W_n ; \; n\geq 0\} \quad {\rm , }\quad \til  \overset{(law)}{=} \sh (\tbar) 
\end{equation}
and that 
\begin{equation}
\label{ttilde2}
\forall u,v\in \til \; , \quad u\leq v \; \Longrightarrow \; \trr_{\Til} (u) \leq \trr_{\Til}(v).
\end{equation}
So, we first obtain $\tau$ by shuffling the GWI-tree $\tbar$ and then by 
identifying the edges of the resulting marked tree that have the same random
marks. We now give
estimates in 
Proposition \ref{ZestimateGW} and in Proposition \ref{factmomentZGWI} on how much
this edge identification does shrink $\til$. Let us introduce some notations:
with any marked forest $F=(f; (m_u, u\in f))$
we associate a collection $(Z_v (F); \, v\in \U)$ of integers defined by 
$$Z_v (F)= \# \{ u\in f: \, \trr_F (u) =v\}.$$ 
Some key estimates in the proof of Theorem \ref{theolim} rely on a precise
computation of the law of the $Z_v
(F)$'s when $F$ is distributed as a GW-forest or a 
GWI-forest. {\bf From now on until the end of the paper} all the GW or
GWI-forests that we consider share the same offspring distribution $\mu
(k)=ud^k$, $k\geq 0$. We set for any $i \in \N^*$ and for any
$x\in [0,1]$
$$ f(x):=\sum_{k\geq 0} ud^k x^k =\frac{u}{1-dx} \quad {\rm and } \quad f_i (x)
:= f(1-a_i +a_i x) .$$
For any $v=m_1\ldots m_n \in \U $ we also define
$$ f_v := f_{m_1} \circ \ldots \circ f_{m_n} \quad {\rm and} \quad a_v := a_{m_1} \ldots a_{m_n} ,$$
with $f_{\varnothing}={\rm Id}$ and $a_{\varnothing }=1$. {\bf We adopt the
  following  convention:} to simplify notation, we do
not distinguish constants in inequalities and we denote them in a generic way 
by a symbol $K_{\alpha , \beta , \ldots }$ meaning that we bound by a positive
constant that only depends on parameters $\alpha , \beta , \ldots $ etc.

  We first 
describe the law of $(Z_v (\f); \, v\in \U)$ with
$\f=(\varphi \, ; (\mu_u , u\in \varphi))$,  
where $\varphi =(\tau_1 , \ldots , \tau_l)$ is a 
forest of $l$ i.i.d  GW($\mu $)-trees and where conditional on $\varphi$ the marks $(\mu_u , u \in
\varphi)$ are taken mutually independent and distributed in accordance with $\a$. 
\begin{proposition}
\label{ZestimateGW} 
\begin{description}
\item{$(i)$} For any $v,w\in \U$,  
$$ \E \left[ x^{Z_{vw}(\f)} \mid Z_{v}(\f) \right]= f_w (x)^{Z_v (\f)} . $$

\item{$(ii)$} Moreover for any $v=m_1 \ldots m_n\in \U$, 
$$1-f_v(1-x)= \frac{x}{ A(v)x+ B(v) }  \quad {\rm with} \quad  1/B(v)= a_v
(d/u)^n $$ 
and 
$$ A(v)= 1+\frac{u}{d} \, \frac{1}{a_{m_1}} + \left( \frac{u}{d}\right)^2 \,\frac{1}{a_{m_1}a_{m_2}} + \ldots +  
\left( \frac{u}{d}\right)^{n-1} \,\frac{1}{a_{m_1} \ldots a_{m_{n-1}}} \; .$$

\item{$(iii)$} For any positive integer $p$, 
$$\E \left[  \sum_{v\in \U} Z_v (\f )^p \right] \leq K_{\a ,p} \, \frac{l^p}{1-d/u } .$$
\end{description}
\end{proposition}

\noindent {\bf Proof :} We first show $(i)$ whose proof reduces to the
``$l=1$'' case 
by an immediate independence argument. Let us take 
$\f=\t_1 =(\tau_1\, ; (\mu_u , u\in \tau_1))$ and  $v\in \U$. Consider the set  $\l_v$ of the
vertices $u\in \tau_1$ satisfying $\trr_{\t_1}(u)=v$. We denote by $u_1<
\ldots <u_{Z_v (\t_1)}$ the elements of $\l_v$ listed in the lexicographical order. 
As a consequence of Remark \ref{indGW}, we see that conditional on $\l_v $ the marked trees 
$(\theta_{u_i} (\t_1)$, $1\leq i\leq Z_v(\t_1))$ are i.i.d. marked trees
distributed as $\t_1$. Observe
next that for any $w\in \U$ 
$$ Z_{vw }(\t_1 ) = \sum_{i=1}^{Z_{v} (\t_1)} \# \{ u\in \theta_{u_i} (\t_1) \, : \; \trr_{\theta_{u_i} (\t_1)} (u) =w \} .$$ 
So we get 
$$\E \left[ x^{Z_{vw (\t_1)}} \mid Z_v (\t_1) \right]= \E \left[  x^{Z_w (\t_1)}  \right]^{ Z_v (\t_1)} .$$
Then it remains to prove: $\E [ x^{Z_w (\t_1) }]= f_w(x)$,  
which follows from iterating the previous identity and from the easy 
observation: $\E [ x^{Z_i (\t_1) }]= f_i (x)$, $i\geq 1$.

The proof of $(ii)$ is a simple recurrence. 
Let us prove $(iii)$: for any positive integer $p$ and any $v=m_1 \ldots m_n \in \U $, we deduce from 
$(ii)$ the following inequality
\begin{eqnarray} 
\label{derivee00}
f_v^{(p)}(1) &=& p\, ! \, \frac{1}{B(v)} \, \left( \frac{A(v)}{B(v)}\right)^{p-1} \\
\label{derivee}  &\leq & p\, ! \,  a_v \left( \frac{d}{u}\right)^{|v|} \left( 1-a_+ \right)^{1-p} ,
\end{eqnarray}
where we have set $a_+ = \max_{i\geq 1} a_i <1$. For any integer $i$ we
denote by $(x)_i$ 
the factorial polynomial $x(x-1)\ldots (x-i+1)$ (with
the convention: $(x)_0=1 $). Check recursively that for any $l, p\geq 1$ and any $h\in C^{\infty}(\R , \R)$,
\begin{equation}
\label{powerderiv}
\frac{\der^p h^l}{\der x^p} = \sum_{j=1}^{p} (l)_j \, h(x)^{l-j} \, Q_{j,p}
(\, h'(x) , \ldots , h^{(p)}(x) \, ) ,
\end{equation}
where the $Q_{j,p}$'s are $j$-homogeneous polynomials with $\N$-valued coefficients that only depend on $j$ and $p$. 
Deduce from (\ref{derivee}) that for any $v\in \U $, 
\begin{eqnarray}
\label{factmomentZ00} 
 \E \left[ (Z_v (\f))_p \right] & =& \frac{\der^p f_v^l}{\der x^p} (1) \\
                             &=& \sum_{j=1}^{p} (l)_j \, Q_{j,p} (\, f_v'
                               (1), \ldots , f_v^{(p)}(1)\, ) \\
                               &\leq & 
                               \sum_{j=1}^{p} (l)_j \, a_v^j \left( \du
                               \right)^{j|v|} \, Q_{j,p} (\, 1! ,
                               2!(1-a_+)^{-1},  \ldots , p!(1-a_+)^{1-p} \, ) \\
\label{factmomentZ}           &\leq & K_{\a,p} \, l^p a_v \left( \du \right)^{|v|}.
\end{eqnarray}
Then, by an easy argument  
\begin{equation}
\E \left[  Z_v (\f)^p \right] \leq  K_{\a , p} \, l^p  a_v \left( \du \right)^{|v|}
\end{equation}
which implies $(iii)$ by the following observation:
\begin{equation}
\label{equila}
\sum_{v\in \U} a_v \left( \du \right)^{|v|} =\sum_{n\geq 0}
\left( \du \right)^{n} \sum_{m_1 , ... ,m_n \geq 1}
a_{m_1} \ldots a_{m_n} = \frac{1}{1-d/u} \; . 
\end{equation} 
\cqfd

\vspace{3mm}

We need similar results for GWI-forests. Let $r$ be some fixed repartition
probability measure. We denote by $\nu$ the corresponding 
immigration distribution given by $\nu (k-1)=\sum_{1\leq j\leq k}
r(k,j)$ , $k\geq 1$. For any $x\in [0,1]$ and any $i\in 
\N^*$ we write 
$$ g(x):=\sum_{k\geq 0} \nu (k) x^k \quad {\rm and } \quad g_i (x)
:= g(1-a_i +a_i x) .$$
Let $\f_0 = (\varphi_0 ; (\mu_u , u \in \varphi_0))$ be a random 
marked GWI-forest whose distribution is characterized as follows: $\varphi_0 =
(\tau_0 , \tau_1, \ldots \tau_l)$, the $\tau_i$'s are mutually independent, 
$\tau_1 , \ldots , \tau_l$ are i.i.d. GW($\mu $)-trees, $\tau_0$ is a GWI($\mu
, r$)-tree and conditional on $\varphi_0$ 
the marks $\mu_u$ are i.i.d. random variables distributed in accordance with
$\a$. For convenience of notation, we set 
$$ u_n^*= u^*_n(\varphi_0) \quad {\rm and }\quad  v^*_n = \trr_{\f_0} (u^*_n)
\; ,\;  n\geq 0 .$$
We also set $\spin =\{ v^*_n i \, , \; i\in \N^* \setminus \{
\mu_{u^*_n}\}\, , n \geq 0\}$ and we define $\s $ as the $\sigma$-field generated 
by the random variables  $(\mu_{u^*_n}\, ; n\geq 0)$ and $(Z_w(\f_0)\, ; w\in \spin )$. 
\begin{proposition}
\label{factmomentZGWI}

$\quad (i)$ Conditional on $\s$, the collection of the $\U$-indexed
processes \\ 
$(\, (Z_{wv} (\f_0) ; \,
v\in \U) \, ; w\in \spin \, )$ are mutually 
independent. Moreover, for any $w\in \spin$, the process $(Z_{wv} (\f_0); v\in \U)$ only 
depends on $\s$ through $Z_{w} (\f_0)$. More precisely, 
$$ (Z_{wv} (\f_0) ; \, v\in \U) \quad {\rm under} \quad \P (\; \cdot \; | w\in \spin \,
;\, Z_{w} (\f_0) =l) \overset{(law)}{=} (Z_v (\f ); \, v\in \U) $$
where $\f=(\varphi \, ; (\mu_u , u\in \varphi))$, where 
$\varphi$ is a sequence of $l$ i.i.d. GW($\mu $)-trees and 
where conditional on $\varphi $ the marks $(\mu_u , u\in \varphi )$ are
i.i.d. distributed in accordance with $\a $. 
\begin{description}
\item{$(ii)$} For any $p\geq 1$, any $n\geq 0$,  
$$\E \left[  Z_{v^*_n } (\f_0)^p \right]  \leq K_{\aa ,p} \, (l+1)^p \max_{0\leq j \leq p} g^{(j)}(1)^p  $$
and for any $i\in \N^*$, 
$$ \E \left[  Z_{v^*_n i} (\f_0)^p \right] \leq K_{\aa ,p} \, a_i \, (l+1)^p \max_{0\leq j \leq p} g^{(j)}(1)^p , $$
\noindent with the convention $g^{(0)}=g$.

\item{$(iii)$} For any $p\geq 1 $ and any $n \geq 0 $,
$$ \E \left[   \sum_{v \in \U} Z_v ([\f_0]_{u^*_n})^p \right] \leq K_{\aa ,
  p} \, \frac{n+1}{1-d/u} (l+1)^p 
\max_{0\leq j \leq p} g^{(j)}(1)^p \,  .$$ 
\end{description}
\end{proposition}

\noindent {\bf Proof :} Set for any $w\in \spin$, $\l_w= \{ u\in \varphi_0 :
\, \trr_{\f_0} (u) =w\}$. Then by definition, $Z_w (\f_0)=\# \l_w $. Check that 
$$ \forall u \neq u'\in \bigcup_{w\in \spin} \l_w \; : \quad  \quad  u , u'\notin \ell_{\infty} (\varphi_0) \quad {\rm and} \quad u\wedge u' \notin \{ u,u'\} .$$ 
These two observations combined with Remark \ref{indGW} imply that 
conditional on $\s $ the marked trees $\theta_u (\f_0)$, $ u \in \bigcup_{w\in \spin} \l_w$ 
are i.i.d. marked GW($\mu $)-trees with independent marks distributed in
accordance with $\a$. This implies $(i)$ thanks to the following equality valid for any $w\in \spin $ and any $v\in \U $:
$$ Z_{wv }(\f_0) = \sum_{u\in \l_w} \# \{ u' \in \theta_u (\f_0 ) \; : \; \trr_{\theta_u (\f_0)} (u')= v\} .$$
Let us prove $(ii)$: Suppose that the word $u^*_n $ is written $l_1 \ldots l_n
\in \U $ for some
nonnegative integers $l_1, \ldots, l_{n}$. 
Consider $u\in \varphi_0$ such that $\trr_{\f_0} (u) =v^*_n$. There
are three
cases: 

$\bullet$ If $|u\wedge u^*_n|=n$ then $u=u^*_n$.

$\bullet$ If $|u\wedge u^*_n|=-1$, then 
$u\wedge u^*_n$ is the fictive root $(-1, \varnothing)$. Thus, the 
ancestor $\varnothing_0$ of the sin-tree $\tau_0$ is not an ancestor of
$u$. It implies 
$$\# \{ u \in \varphi_0 \; :\; |u\wedge u^*_n|=-1 \quad {\rm and} \quad
\trr_{\f_0} (u) =v^*_n \} =\sum_{j=1}^l  Z_{v^*_n} \left(
  \theta_{\varnothing_j} (\f_0) \right) .$$

$\bullet$ If $u\wedge u^*_n=u^*_k$ with $0\leq k<n$, we can find some 
$j\in \{ 1, \ldots , k_{u^*_k} (\varphi_0) \}$ with $j\neq l_{k+1}$ and some $u'\in \U $ such that 
$$ u=u^*_k j u' \quad , \quad \mu_{u^*_k j}=\mu_{u^*_{k+1}} \quad {\rm and}
\quad \trr_{\theta_{u^*_k j}(\f_0)} (u') =w^*_{k+1} , $$
where $w^*_{k+1} \in \U$ stands for the word $\mu_{u^*_{k+2}} \ldots
\mu_{u^*_n} \in \U $, with the convention $w^*_{n}=\varnothing $.

\vspace{3mm}

Now, set for any $0\leq k <n$,  
$$ E_{k+1} =\{  j\in \{ 1, \ldots , k_{u^*_k} (\varphi_0) \} \; : \; j\neq
l_{k+1} \quad {\rm and} \quad \mu_{u^*_{k} j}= \mu_{u^*_{k+1}} \} .$$
The combination of the three preceding cases implies that
$$ Z_{v^*_n} (\f_0)=  1 + \sum_{j=1}^l  Z_{v^*_n} \left(
  \theta_{\varnothing_j} (\f_0) \right) + \sum_{k=0}^{n-1}  
\sum_{j\in E_{k+1}} 
\# \{ u' \in \theta_{u^*_{k} j } (\varphi_0) \; : \; \trr_{\theta_{u^*_{k} j } (\f_0)} (u') =w^*_{k+1} \} .$$

Set $\kappa_k = \# E_k $ , $1\leq k\leq n$, 
$\kappa_0=l$ and $w^*_0=v^*_n$. Then, by $(i)$ we get 
\begin{eqnarray}
\label{spindecZ0}
\E \left[ r^{Z_{v^*_n} (\f_0)} 
\mid (\kappa_{k+1},\mu_{u^*_{k+1}} )_{0\leq k\leq n-1}
\right] & = & r f_{v^*_n} (r)^l \prod_{k=0}^{n-1}
f_{w^*_{k+1}} (r)^{\kappa_{k+1} } \\
\label{spindecZ} &= &  r\prod_{k=0}^n
f_{w^*_k} (r)^{\kappa_k } . 
\end{eqnarray}
It also follows from the 
previous observations that $\kappa_1 , \ldots ,\kappa_n$ are mutually
independent with the same distribution specified by 
\begin{equation}
\label{loikappa}
\E\left[ x^{\kappa_1}\right] = \E\left[ g_{\mu_{u^*_0}} (x)\right]= 
\sum_{i\in  \N^*} a_i g_{a_i} (x) .
\end{equation}
From Proposition
\ref{ZestimateGW} $(ii)$ we get a.s.
$$ f_{w^*_k} (1+z) = 1+ \frac{z}{B(w^*_k)} \, \left( 1
  -\frac{A(w^*_k)}{B(w^*_k)} \, z\right)^{-1} $$
and since  
$$\frac{A(w^*_k)}{B(w^*_k)} \leq (\du a_+)^{|w^*_k|} + 
\ldots + \du a_+ \leq  (1-a_+)^{-1} , $$
$f_{w^*_k} (1+z)$ has a.s. a power serie expansion with a 
radius of convergence greater than $1-a_+>0$. Then, for any $|z| < 1-a_+$
we can write 
$$ f_{w^*_k} (1+z)^{\kappa_k } = 1+ \sum_{p\geq 1} D^{(k)}_p z^p  \quad 
{\rm with}
\quad  D^{(k)}_p =
\frac{1}{p!} \frac{\der^p f_{w^*_k}^{\kappa_k }}{\der z^p} (1).$$
Deduce from (\ref{factmomentZ})
\begin{eqnarray}
\label{estimateDp0}
0 \leq \frac{1}{p!} \frac{\der^p f_{w^*_k}^{\kappa_k }}{\der z^p} (1) &\leq & 
K_{p, \aa } \, \kappa_k^p \, a_{w^*_k} \, \left( \du \right)^{|w^*_k|} \\
\label{estimateDp} &\leq &  K_{p, \aa } \, \kappa_k^p \,  a_+^{n-k} .
\end{eqnarray}
Then observe that 
$$  \prod_{k=0}^n f_{w^*_k} (1+z)^{\kappa_k} = 1 +\sum_{p\geq 1} D_p z^p \quad
, \; |z| < 1-a_+$$
where 
$$ D_p = \sum_{\p \subset \{ 0, \ldots, n \} } \; \; 
\sum_{\substack{ \sum_{k\in \p} q_k =p \\ q_k \geq 1 }} \; \;  \prod_{k\in \p} D_{q_k}^{(k)} .$$
Set $D_0 =1$ and deduce from (\ref{spindecZ}) 
\begin{equation}
\label{inegzero}
\E \left[ \left( Z_{v^*_n} (\f_0) \right)_p \mid (\kappa_k,\mu_{u^*_{k}} )_{0\leq k\leq n} \right] = p! \, \left( D_p + D_{p-1} \right)
\end{equation}
Use (\ref{estimateDp}) and the independence of the $\kappa_i$'s to get 
$$\E \left[ D_p \right] \leq 
\sum_{\p \subset \{ 0, \ldots, n \} } 
\sum_{\substack{ \sum_{k\in \p} q_k =p \\ q_k \geq 1 }} \; \prod_{k\in \p}   
K_{q_k, \aa } \E \left[ \kappa_k^{q_k} \right] \, a_+^{n-k}. $$
If $\p\subset\{ 0, \ldots, n \} $ and $ \sum_{k\in \p} q_k =p $ with $q_k
\geq 1 $, $k \in \p$, then $\# \p \leq p $ and $q_k \leq p $ for any 
$k\in \p $. Thus, 
$$\prod_{k\in \p}   \E \left[ \kappa_k^{q_k} \right] \leq (l+1)^p  \left( 1\vee \E [ \kappa_1^p ]\right)^p $$
since $\kappa_1 , \ldots , \kappa_n $ are identically distributed and $\kappa_0 =l $. Deduce
from (\ref{loikappa}): 
$$ 1 \vee \E [\kappa_1^p ] \leq K_{\a,p} \max \left( 1, g' (1) , \ldots , g^{(p)}(1) \right) .$$
Thus, 
\begin{equation}
\label{ineg}
\E \left[ D_p \right] \leq 
 K_{\a , p} (l+1)^p \max_{0 \leq j\leq p} g^{(j)}(1)^p  
\end{equation}
since 
$$ 
\sum_{\p \subset \{ 0, \ldots, n \} } 
\sum_{\substack{ \sum_{k\in \p} q_k =p \\ q_k \geq 1 }} \; \prod_{k\in \p}
 a_+^{n-k} \leq  K_{p} \, (1-a_+)^{-p} .$$
Then by (\ref{inegzero}) and an easy argument 
\begin{equation} 
\label{momentdeux}
\E \left[  Z_{v^*_n } (\f_0)^p \right] \leq K_{\aa ,p} \, (l+1)^p \max_{0\leq j \leq p} g^{(j)}(1)^p .
\end{equation}

To achieve the proof
of $(ii)$, we set $\l_{v^*_n}=\{ u\in \varphi_0 \setminus\{ u^*_n\} \; :\;
\trr_{\f_0} (u) = v^*_n\}$. Recall that any $u\in \l_{v^*_n}$ has offspring distribution $\mu$ and that $u^*_n$ 
has offspring distribution $\nu$. Since $Z_{v^*_{n} } =1+ \# \l_{v^*_n}$, we get for any $i\in \N^*$,
$$ \E \left[ x^{Z_{v^*_{n}i }} \mid \l_{v^*_n}, \mu_{u_{n+1}^*} \right]= x^{\un_{\{ \mu_{u_{n+1}^* }=i \}}} 
f_i(x)^{Z_{v^*_{n} } -1} g_i(x) . $$
Thus, 
$$ \E \left[ x^{Z_{v^*_{n}i }} \mid  Z_{v^*_{n}}=k +1\right]= (1-a_i +a_i x) f_i(x)^{k} g_i(x) . $$
By differentiating $p$ times at $x=1$ we get  
$$ \E \left[ \left( Z_{v^*_{ni} } \right)_p \mid  Z_{v^*_{n}}=k +1\right]= \frac{{\rm d}^p f_i^k g_i}{{\rm d} x^p}(1)+ p a_i
\frac{{\rm d}^{p-1} f_i^k g_i}{{\rm d} x^{p-1}} (1) $$
Now observe that for any $q\geq 0$,  $g_i^{(q)}(1)= a_i^{q} g^{(q)}(1)$ and 
$$ \frac{{\rm d}^q f_i^k }{{\rm d} x^q}(1)= (a_i d/u)^q (k+q -1)_q ,$$
by a simple computation. Thus, 
\begin{eqnarray*}
 \frac{{\rm d}^p f_i^k g_i}{{\rm d} x^p}(1) &=&  \sum_{q=0}^p \frac{p!}{q!(p-q)!}a_i^{p-q} 
g^{(p-q)}(1) (a_i d/u)^q (k+q -1)_q     \\
 &\leq & K_p a_i^p (k+p)_p  \max_{1\leq j \leq p} g^{(j)}(1) .
\end{eqnarray*}
Consequently, 
$$\E \left[ \left( Z_{v^*_{ni} } \right)_p \mid  Z_{v^*_{n}}\right] \leq  K_p a_i^p ( Z_{v^*_{n}} - 1+p)_p  \max_{1\leq j \leq p} 
g^{(j)}(1) .$$
Since $p\geq 1$ and by (\ref{momentdeux}) we get 
$$\E \left[ \left( Z_{v^*_{ni} } \right)_p \right] \leq  K_p a_i  (l+1)^p  \max_{1\leq j \leq p} 
g^{(j)}(1)^p $$
which easily implies the second inequality of Lemma \ref{factmomentZGWI} ($ii$).

\vspace{3mm}

We now prove $(iii)$: first observe that the decomposition 
$$ \sum_{v \in \U} Z_v ([\f_0]_{u^*_n})^p = e_1 + e_2 + e_3 $$
holds with 
$$ e_1 = \sum_{\substack{ w\in \spin \\  |w| \leq n} } \sum_{v\in \U} Z_{wv}([\f_0]_{u^*_n})^p ,$$
 
$$ e_2 = \sum_{0\leq k < n} Z_{v^*_k} ([\f_0]_{u^*_n})^p \quad {\rm and } \quad 
e_3 =\sum_{v\in \U}   Z_{v^*_n v} ([\f_0]_{u^*_n})^p .$$
Note for any $v\in \U $ and for any $w\in \spin $ such that $|w| \leq n $ that
$Z_{wv}([\f_0]_{u^*_n})= Z_{wv}(\f_0)$. Then by Proposition
\ref{ZestimateGW} $(i)$ and Proposition \ref{factmomentZGWI} $(iii)$ 
\begin{equation}
\label{e1}
\E \left[ e_1 \mid \s \right] \leq K_{\aa , p} \sum_{\substack{ w\in \spin \\
    |w| \leq n} } \frac{Z_w (\f_0)^p}{1-d/u} \leq   K_{\aa , p} (1-d/u)^{-1} 
\sum_{k=0}^{n}
\sum_{i\in \N^*} Z_{v^*_ki} (\f_0 )^p .
\end{equation}
We then deduce from the second inequality of Lemma \ref{factmomentZGWI}
($ii$)
\begin{equation}
\label{e1bis}
\E \left[ e_1 \right] \leq K_{\aa , p} \,  (l+1)^p \max_{0\leq j \leq p} g^{(j)}(1)^p \frac{n+1}{1-d/u}.
\end{equation}
Observe next that $Z_{v^*_k }([\f_0]_{u^*_n})= Z_{v^*_k}(\f_0)$. Then by
the first inequality of $(ii)$, we get 
\begin{equation}
\label{e2}
\E \left[ e_2 \right] \leq  K_{\aa , p} \, n \; (l+1)^p \max_{0\leq j \leq p} g^{(j)}(1)^p .
\end{equation}
To bound $\E \left[ e_3 \right]$, note that conditional on $Z_{v^*_n} ([\f_0]_{u^*_n})=l$, the
process $(Z_{v^*_n v} ([\f_0]_{u^*_n}) \, ; v\in \U )$
is distributed as $(Z_v (\f ); \, v\in \U)$ where 
$\f=(\varphi \, ; (\mu_u , u\in \varphi))$, 
where $\varphi$ is a sequence of $l$ independent GW($\mu $)-trees and 
where conditional on $\varphi $ the marks $(\mu_u , u\in \varphi )$ are
i.i.d. random variables 
distributed in accordance with $\a $. Thus, by Proposition 
\ref{ZestimateGW} : 
$$ \E \left[ e_3  \mid Z_{v^*_n} ([\f_0]_{u^*_n}) \right] \leq K_{\aa, p} \frac{ Z_{v^*_n} ([\f_0]_{u^*_n})^p}{1-d/u} .$$
Now observe that $Z_{v^*_n} (\f_0 )=Z_{v^*_n} ([\f_0]_{u^*_n})$ and use
(\ref{momentdeux}) to get 
\begin{equation}
\label{e3}
\E \left[ e_3 \right] \leq  K_{\aa , p} \, 
\max_{0\leq j \leq p} g^{(j)}(1)^p  \, \frac{(l+1)^p}{1-d/u}.
\end{equation}
Then, $(iv)$ follows by adding (\ref{e1bis}), (\ref{e2}) and (\ref{e3}). \cqfd

\subsection{Proof of Theorem \ref{theolim}. }
Let us first explain why Theorem \ref{theolim} reduces to
a convergence for finite trees: we restore $\ee $ in the random variables: 
$\Tbaree=(\tbaree \, ; (\mubar_u , u \in \tbaree))$ and $\Tilee=\sh( \Tbaree ) =(\tilee\, ; (\muti_u , u \in \tilee))$.
 For any positive real 
number $x$ we set 
$\xxe =\lfloor x / \ee \rfloor $ and we define $\ztileex= \sup \{ n\geq 0 :
\, |W_n^{\ee}| \leq \xxe \}$. 
As explained in the introduction, we associate a unique 
finite ordered rooted tree $\txe$ with the
subtree $\{ W_n^{\ee } ; \, 0\leq n \leq \ztileex
\} \subset \U$. Observe that in general $\txe \neq [\tte ]_{u^*_{\xxe } (\tte )}$,
however $\txe$ and $\tte $ coincide up to level $\xxe$: 
\begin{equation}
\label{coincide}
[\txe]_{\xxe} = [\tte]_{\xxe} .
\end{equation}
The following proposition asserts that the convergence of $\tte$ is equivalent
to the convergence of the $\txe $'s for all $x>0$. For convenience of notation,   
we set $\zxe = 2\ee^2 \# \txe $ and 
$$ H_s (x, \ee)=  \ee \,\un_{ [0 , \zxe) } (s) \, H_{\lfloor s/ 2\ee^2 \rfloor }( \txe)
\quad {\rm and }\quad C_s(x, \ee) =\ee \, \un_{ [0, \zxe -2\ee^2 ] } (s)\, 
C_{s/\ee^2  }( \txe ) . $$
We also define the limiting process by  
$$ D^{(x)}_s = \un_{[0, \sigma_x ]} (s)
D_s + \un_{[\sigma_x , \infty )} (s) D^{\bullet}_{(\zeta_x -s)_+} ,$$ 
where $\zeta_x =\sigma_x + \sigma^{\bullet}_x $ with $\sigma_x $ 
(resp. $\sigma^{\bullet }_x$) $\, =\sup \{ s\geq 0: \, D_s \, ({\rm resp.} \; D^{\bullet}_s ) \leq x \}$.
\begin{proposition}
\label{theolimfini}
Theorem \ref{theolim} is implied by any of the following equivalent convergences 
 
$$ (i) \quad \quad \forall x > 0 \; , \quad  C(x, \ee ) \xrightarrow[\ee\rightarrow 0]{\quad ({\rm d}) \quad } \left(  2D^{(x)}_{\gamma s}; \; s\geq 0 \right) $$
$$ (ii) \quad \quad \forall x > 0 \; , \quad  H(x, \ee ) \xrightarrow[\ee\rightarrow 0]{\quad ({\rm d}) \quad } \left( 2 D^{(x)}_{\gamma s}; \; s\geq 0 \right) .$$
\end{proposition}

\noindent {\bf Proof:} The proof of $(ii) \Longrightarrow (i)$ can be copied
from the proof of Theorem 2.4.1 \cite{DuLG}. It relies on 
formula (\ref{contourvsheight}) that makes the contour process of a finite
ordered rooted tree an explicit functional of the corresponding height 
process. 
Since  (\ref{contourvsheight}) also holds for contour processes of
sin-trees, similar arguments work to show that Theorem \ref{theolim} $(ii)$
implies Theorem \ref{theolim} $(i)$. Let us prove that
Proposition \ref{theolimfini} $(i)$ implies Proposition \ref{theolimfini} $(ii)$: 
Recall from (\ref{contourvsheight}) that 
$$ H_n (\txe) = C_{2n-H_n (\txe )} (\txe ). $$
So, if we denote by $S(\ee)$ the maximal height of $\txe $ we get 
$$ \sup_{n< \# \txe } \mid H_n (\txe ) -C_{2n} (\txe ) \mid \; \leq
\max_{|n-n'| \leq S(\ee)} \mid C_n (\txe ) - C_{n'} (\txe) \mid ,$$
which implies after scaling
$$ \sup_{s \leq \zxe } \mid H_s (x, \ee ) -C_{s} (x, \ee ) \mid \; \leq \max_{|s-s'| \leq \ee^2 S(\ee)} \mid C_s (x, \ee ) - C_{s'} (x, \ee) \mid . $$
Proposition \ref{theolimfini} $(i)$ implies that $\ee S(\ee)$ converges in distribution to the
supremum of $D^{(x)}$ that is a.s. finite. Thus, the right member of the latter inequality converges to zero in
probability and Proposition \ref{theolimfini} $(ii)$ follows. 
A similar argument show that Theorem \ref{theolim} $(i)$ implies 
Theorem \ref{theolim} $(ii)$. 
Now, the proof will be achieved if we show that Proposition \ref{theolimfini} $(ii)$
implies Theorem \ref{theolim} $(ii)$: Assume that Proposition \ref{theolimfini}
$(ii)$ is true and deduce from (\ref{coincide}) that
\begin{equation}
\label{hittingtime} 
\quad \left( \ee H_{\lfloor s\wedge {\bf e}_{x, \ee} /2\ee^2 \rfloor }(\tau_{\ee})
  \; , \; \ee H^{\bullet}_{\lfloor 
s\wedge {\bf e}^{\bullet}_{x, \ee} /2\ee^2 \rfloor} (\tau_{\ee}) 
\right)_{s\geq 0}
\xrightarrow[\ee\rightarrow 0]{\quad ({\rm d}) \quad} \left( 2D_{\gamma
  (s\wedge {\bf e}_x)}\; , \; 2D^{\bullet}_{\gamma (s\wedge {\bf e}^{\bullet}_x) }\right)_{s\geq 0}, 
\end{equation}
where ${\bf e}_{x, \ee}= \inf \{ n \geq 0 : \, H_n (\txe ) \geq \xxe \} $ and ${\bf e}_x =\inf \{ s \geq 0 : \, D_s \geq x \} $
with similar definitions for ${\bf e}^{\bullet}_{x, \ee}$ and 
${\bf e}^{\bullet}_x$. Observe that Proposition \ref{theolimfini} $(ii)$ implies for any
$x>0$ that $({\bf e}_{x, \ee}, {\bf e}^{\bullet}_{x, \ee})$ converges in
distribution to $({\bf e}_{x}, {\bf e}^{\bullet}_{x})$. Since 
${\bf e}_{x}$ and ${\bf e}^{\bullet}_{x}$ a.s. go to infinity with $x$, we
then get for any $M>0$,
$$ \lim_{x \rightarrow \infty } \limsup_{\ee \rightarrow 0} 
\P \left( e_{x, \ee} \leq M \, ;  \, e^{\bullet}_{x, \ee} \leq M \right) =0 ,$$
which implies Theorem \ref{theolim} $(ii)$ by (\ref{hittingtime}) and by standard arguments. \cqfd

\vspace{3mm}


We define $ \tbareex = [\tbaree]_{u^*_{\xxe} (\tbaree)} $
and $\tileex =[\tilee]_{u^*_{\xxe} (\tilee)}$ and we also set
$$ \Tbareex= [\Tbaree]_{u^*_{\xxe} (\tbaree)} = \left( \tbareex \, ; (\mubar_{u} , u \in \tbareex ) \right)  
\quad {\rm and } \quad \Tileex =[\Tilee]_{u^*_{\xxe} (\tilee)}= \left( \tileex \, ; (\muti_{u} , u \in \tileex ) \right) .$$ 
By definition, $\# \tileex = \# \tbareex =\ztileex$. Deduce from (\ref{ttilde1}) and (\ref{ttilde2}) 
\begin{equation}
\label{combifond1}
\trr_{\Tilee}(\tileex)=\{ W_n^{\ee} ; \; 0\leq n \leq \ztileex \} \quad {\rm , }\quad \tileex  \overset{(law)}{=} \sh (\tbareex) 
\end{equation}
and 
\begin{equation}
\label{combifond2}
\forall u,v\in \tileex \quad  : \quad u\leq v \; \Longrightarrow \; \trr_{\Tilee} (u) \leq \trr_{\Tilee}(v) .
\end{equation}
By Proposition \ref{theolimfini}, Theorem \ref{theolim} reduces to prove that for any $x > 0$:
\begin{equation}
\label{mainconv}
H(x, \ee ) \xrightarrow[\ee\rightarrow 0]{\quad ({\rm d}) \quad } \left( 2 D^{(x)}_{\gamma s}; \; s\geq 0 \right).
\end{equation}
The first step of the proof of (\ref{mainconv}) is a limit theorem for
$\tileex$: let us set for any $s\in [0, \infty )$
$$  \tilH_s (x, \ee)=  \ee \,\un_{ [0 \, , \, 2\ee^2 \# \tileex ) } (s) \, H_{\lfloor s/ 2\ee^2 \rfloor }( \tileex)
\quad {\rm and }\quad \tilC_s(x, \ee) =\ee \, \un_{ [0\, , \, 2\ee^2(\# \tileex -1) ] } (s)\, 
C_{s/\ee^2  }( \tileex ).$$
\begin{lemma}
\label{convtileex}
$$ (i) \quad \quad \forall x > 0 \; , \quad  \tilC (x, \ee ) \xrightarrow[\ee\rightarrow 0]{\quad ({\rm d}) \quad } \left(  2D^{(x)}_{s}; \; s\geq 0 \right) $$
$$ (ii) \quad \quad \forall x > 0 \; , \quad  \tilH (x, \ee ) \xrightarrow[\ee\rightarrow 0]{\quad ({\rm d}) \quad } \left( 2 D^{(x)}_{s}; \; s\geq 0 \right) $$
\end{lemma}
\noindent {\bf Proof :} Deduce from (\ref{spinaldec3})
\begin{equation}
\label{spinaltilee}
H_n (\tilee) = n - {\bf p}_{\ee} (n) + H_{{\bf p }_{\ee} (n)} ( f (\tilee)) \;
, \quad n\geq 0 .
\end{equation}
Recall that $f(\tilee )$ stands for the forest composed by the bushes rooted
at the left hand of the infinite line of descent of $\tilee$ and 
that ${\bf p}_{\ee} (n) $ is given by ${\bf p}_{\ee}(n)= \inf \{ p\geq 0\; :\; {\bf n}_{\ee}(p)\geq n \}$ where
$$ {\bf n}_{\ee}(p)=p+\inf \{ k\geq 0 :\; L_k (\tilee)> -\inf_{j\leq p } W_j ( f(\tilee) )\} $$
with $L_n (\tilee)=(l_1 (\tilee)-1) + \ldots + (l_n(\tilee)-1)$ , $n\geq 1$ and $L_0 (\tilee)=0$. By Remark \ref{shuffleGWI}, $\tilee$ is a 
GWI($\mu , r'$)-tree with $r'(k,l)= ud^{k-1} /k$, $1\leq l \leq k$. Thus by Remark \ref{spinalGWI}:
\begin{description}
\item{$-$} the two forests $(f(\tilee),  f(\tilee^{\bullet}))$ are independent of $(L(\tilee) , L(\tilee^{\bullet}))$; 
\item{$-$} $f(\tilee)$ and $ f(\tilee^{\bullet})$ are two mutually independent
  sequence of i.i.d GW($\mu $)-trees;
\item{$-$} $(L(\tilee) , L(\tilee^{\bullet}))$ is a $\N \times \N $-valued
  random walk whose jump distribution is given by
$$ \P \left( L_{n+1} (\tilee) -L_n (\tilee)= l \, ; \,  L_{n+1} (\tilee^{\bullet}) -L_n (\tilee^{\bullet})=l' \right)= \frac{1}{l +l'+1}\, ud^{l+ l'}.$$
\end{description}
Check first that $\E [L_n (\tilee)]= \E [L_n (\tilee^{\bullet})]=nd/2u$, which implies 
\begin{equation}
\label{spintaux}
\left( \ee L_{s/\ee}(\tilee) \; , \; \ee L_{s/\ee} (\tilee^{\bullet}) \right)_{s\geq 0}
\xrightarrow[\ee\rightarrow 0]{\quad ({\rm d}) \quad} \left( s/2 \; , \;  s/2\right)_{s\geq 0} .
\end{equation}
Next, we need to prove the joint convergence of $(\ee \, H_{\lfloor
  s/2\ee^2\rfloor}((\tilee)), \ee \, 
V_{\lfloor s/2\ee^2 \rfloor}(f(\tilee)))$: We know 
from Remark \ref{harrisrem} that $(V_p (f(\tilee )) ; \; p\geq 0)$ is a random walk with jump distribution given by $\rho(k)=ud^{k+1}$, $k\geq -1 $.  
An elementary computation implies for any $\lambda \in \R$ that 
$$ \E \left[ \exp \left( i\lambda  \ee \, V_{\lfloor s/2\ee^2 \rfloor}(f(\tilee)) \right)\right]= 
\exp \left( -\frac{s\lambda^2}{2} -2i\lambda s \right) + o(1) $$
and by standard arguments   
\begin{equation} 
\label{convwalk}
 \left( \ee \, V_{\lfloor s/2\ee^2 \rfloor}(f(\tilee)) \right)_{s\geq 0}
\xrightarrow[\ee\rightarrow 0]{\quad ({\rm d}) \quad} B^{(-2)} 
\end{equation}
(see for instance Theorem 2.7 \cite{Sko}). We then use Theorem 2.3.1
\cite{DuLG} that asserts that under 
(\ref{convwalk}) the following joint convergence 
\begin{equation} 
\label{conjointe}
 \left( \ee \,  H_{\lfloor s/2\ee^2\rfloor}(f(\tilee))\; ,\; \ee \,  V_{\lfloor s/2\ee^2 \rfloor}(f(\tilee)) \right)_{s\geq 0}
\xrightarrow[\ee\rightarrow 0]{\quad ({\rm d}) \quad} \left( 2(B^{(-2)} -I^{(-2)}) \; , \; B^{(-2)} \right)
\end{equation}
holds provided that for any $\delta >0 $, 
\begin{equation}
\label{hypo}
 \liminf_{\ee \rightarrow 0} \left( f_{\lfloor \delta/\ee \rfloor} (0) \right)^{\lfloor 1/\ee \rfloor} >0 
\end{equation}
(recall that $f_n $ is recursively defined by $f_n =f_{n-1} \circ f $). Check that 
$$ f_n (x) = \frac{u}{d} \, \frac{1-\left(u/d \right)^{n} -x (1-\left(u/d \right)^{n-1} )}{1-\left( u/d \right)^{n+1} -x 
\left( 1-\left( u/d \right)^{n} \right)} .$$
Then, 
$$  \lim_{\ee \rightarrow 0} 
\left( f_{\lfloor \delta/\ee \rfloor} (0) \right)^{\lfloor 1/\ee \rfloor} =
\exp \left( -\frac{4}{e^{4\delta}-1} \right) \; >\; 0 $$ 
and (\ref{conjointe}) follows from (\ref{hypo}). Recall notation $\alpha$
from Section 2.2 and observe that   
$$\ee \aal (\lfloor s/2 \ee^2 \rfloor )= 
\inf \{ s'\geq 0 :\; \ee  L_{\lfloor s'/ \ee \rfloor} (\tilee)> -\inf_{r\leq \lfloor s/ 2\ee^2 \rfloor } 
\ee \, V_{\lfloor r/ 2\ee^2 \rfloor} ( f(\tilee) ) \} .$$ 
Deduce from (\ref{spinaldec1}) that:
$$ 2\ee^2{\bf n}_{\ee}( \lfloor s/ 2\ee^2 \rfloor)=  2\ee^2 \lfloor s/
  2\ee^2 \rfloor  +  2\ee^2  \aal (\lfloor s/2 \ee^2 \rfloor ) .$$ 
Then by (\ref{spintaux}) and (\ref{conjointe})
$$\left( \ee \aal (\lfloor s/2 \ee^2 \rfloor ) \right)_{s\geq 0}  
\xrightarrow[\ee\rightarrow 0]{\quad ({\rm d}) \quad} -2I^{(-2)}
\quad {\rm and} \quad (2\ee^2 {\bf n}_{\ee} (\lfloor s/ 2\ee \rfloor)\, ; s\geq 0) \overset{{\rm (d)}}{\rightarrow} (s\, ; s\geq 0)
.$$
Thus, $(2\ee^2 {\bf p}_{\ee} (\lfloor s/ 2\ee
\rfloor)\, ; s\geq 0) \overset{{\rm (d)}}{\rightarrow} (s\, ; s\geq 0)
$ and (\ref{spinaldec4}) combined with the convergence of 
${\bf p}_{\ee}$ 
and (\ref{spinaltilee}) imply 
$$ \left( \ee H_{\lfloor s/ 2\ee^2 \rfloor} (\tilee )\right)_{s \geq 0} \xrightarrow[\ee\rightarrow 0]{\quad ({\rm d}) \quad} 
\left( 2B_s^{(-2)}-4I_s^{(-2)} \right)_{s\geq 0} \;= \; 2D .$$
The joint convergence (\ref{spintaux}) combined with the independence of $ f (\tilee)$ and $ f (\tilee^{\bullet })$ also implies 
$$ \left( \ee H_{\lfloor \cdot / 2\ee^2 \rfloor} (\tilee ) 
\; , \; \ee H_{\lfloor \cdot / 2\ee^2 \rfloor} (\tilee^{\bullet} ) \right)_{s \geq 0}  
 \xrightarrow[\ee\rightarrow 0]{\quad ({\rm d}) \quad} \left( 2D \; , \; 2D^{\bullet} \right) .$$
Use (\ref{contourvsheight}) and arguments similar to those used 
in the proof of Theorem 2.4.1 \cite{DuLG} to get  
\begin{equation}
\label{convjointcont}
\left( \ee \, C_{ s/ \ee^2 } (\tilee ) \; , \; \ee \, C^{\bullet }_{ s/ \ee^2 } (\tilee ) \right)_{s \geq 0} 
\xrightarrow[\ee\rightarrow 0]{\quad ({\rm d}) \quad} 
\left( 2D \; , \; 2D^{\bullet} \right).
\end{equation}
Set $\widetilde{\sigma }_{x , \ee}= 
\sup \{ s\geq 0 \, : \; C_s (\tilee ) \leq \xxe \}$ and define 
$ \widetilde{\sigma }_{x , \ee}^{\bullet}$ in a similar way. 
Recall notations $\sigma_x$, $\sigma_x^{\bullet}$ and $D^{(x)}$ introduced before Proposition \ref{theolimfini} and deduce from 
(\ref{convjointcont}) that 
$$ \left( \widetilde{\sigma }_{x ,\ee} 
\; , \; \widetilde{\sigma }_{x , \ee}^{\bullet} \right) 
\xrightarrow[\ee\rightarrow 0]{\quad ({\rm d}) \quad} \left( \sigma_x \; , \; \sigma_x^{\bullet} \right) .$$
It easily implies Lemma \ref{convtileex} $(i)$ by (\ref{cutcontour2}). Then,
argue exactly as in the proof of Proposition \ref{theolimfini} to 
deduce Lemma \ref{convtileex} $(ii)$ from Lemma \ref{convtileex} $(i)$. \cqfd

\vspace{6mm}

We now have to prove Proposition \ref{theolimfini} $(ii)$. In one
part of the proof we adapt Aldous' approach (Theorem 20 \cite{Al2}) and we
get estimates for the tree 
$\txe $ reduced at certain random times. The main technical difficulty is Lemma \ref{convproba} 
that asserts that these random times are asymptotically uniformly distributed.
Let us first define these random times: Let $(\u_i \, ; i\geq 1 )$ be a sequence of i.i.d. random variables independent of 
$W^{\ee}$ and uniformly distributed on $(0,1)$. Let $u_0 = \varnothing < u_1 < \ldots < u_{\# \tileex -1} $ be the vertices of $\tileex $ listed in the 
lexicographical order. We set 
$$ U_i (x,\ee)= u_{\lfloor \u_i \# \tileex \rfloor } \quad {\rm and } \quad
V_i (x, \ee) = 
\trr_{\Tileex} \left( U_i(x, \ee) \right) \in \U .$$
Then $V_i (x, \ee)\in  \{ W^{\ee}_n \, ; 0\leq n \leq \ztileex \} $ and the row of the corresponding vertex in $\txe $ is given by 
$$ \Vbar_i (x,\ee)= \sum_{\substack{v\in \U  \\ v\leq V_i(x, \ee)}} \un_{\{ Z_v (\Tileex ) >0 \}} . $$
The key argument is the following Lemma that is proved in the next section.
\begin{lemma}
\label{convproba}
For any $i\geq 1$, the following convergence holds in probability: 
$$ \ee^2 \left( \Vbar_i (x, \ee) - \frac{1}{\gamma} \, \u_i \, \# \, \tileex \right) \xrightarrow[\ee\rightarrow 0]{\quad  \quad} 0. $$ 
\end{lemma}

\noindent
From now on until the end of the section we assume that Lemma \ref{convproba} is 
true and we prove Proposition \ref{theolimfini} $(ii)$:
Fix $x>0$ and set for any $\delta > 0$:
$$ \omega \left( H(x, \ee) , \delta \right) = \sup \{ \left|  H_s (x,\ee) -
  H_{s'} (x, \ee) \right| \; ; \; \left| s-s' \right| \leq \delta \} .$$
We first prove tightness for $H(x, \ee)$,
$\ee >0$: By a standard criterion (see for instance Corollary 3.7.4 \cite{EK})
we only need to prove 
$$ {\rm(T1) } \quad \quad \lim_{M \rightarrow \infty} \liminf_{\ee \rightarrow
  0} 
\P \left( \sup_{s\geq 0} H_s (x, \ee ) \leq M \right) \, = \, 1 $$
and for any $\eta > 0$,  
$$ {\rm (T2)} \quad \quad \lim_{\delta \rightarrow 0} \limsup_{\ee \rightarrow 0} \P \left( \omega \left( H(x, \ee) , \delta \right) > \eta \right) \, = \, 0. $$ 
{\bf Proof of (T1) :} Note that the mapping $\trr$ preserves height. So, 
we get 
\begin{equation*}
\sup_{s\geq 0} H_s (x, \ee) = \ee \sup \{  \trr_{\Tileex} (u) \; :\;
u\in \tileex \} =  \sup_{0 \leq s < 2\ee^2 \# \tileex} \ee H_{\lfloor s/ 2\ee^2 \rfloor }
 (\tileex ) 
\end{equation*} 
which is a tight family of random variables by Lemma \ref{convtileex}. \cqfd 

\vspace{3mm}

\noindent {\bf Proof of (T2) :} Let $k$ be a positive integer and $p$ be a permutation of $\{1, \ldots, k \}$ such that 
$V_{p(1)} (x, \ee) \leq \ldots \leq V_{p(k)} (x, \ee) $ in $\U $. It implies 
$$ \Vbar_{p(0)} (x, \ee ) \leq \Vbar_{p(1)} (x, \ee) \leq \ldots \leq \Vbar_{ p(k)} (x,\ee) \leq \Vbar_{p(k+1)} (x, \ee) $$
where we set $0=\Vbar_{p(0)} (x, \ee)$ and $\# \txe =
\Vbar_{p(k+1)} (x, \ee) $. 
We first need to get an upperbound for the quantities $q_i$ defined
for any $0\leq i \leq k $ by 
$$ q_i = \sup \left\{ \, \left|  H_n (\txe ) - H_{\Vbar_{p(i)} (x, \ee)}  (\txe ) \right| \; ; \;  \Vbar_{p(i)} (x, \ee)\leq n \leq \Vbar_{p(i+1)} (x, \ee)
\right\} .$$
Observe that $q_i $ can be rewritten  
\begin{equation}
\label{rewriteq}
q_i = \sup \left\{ \left|  \, |v| - \left| V_{p(i)} (x, \ee) \right| \, \right| \; ;\; v \in \trr_{\Tileex} (\tileex) \quad {\rm and } \quad 
V_{p(i)} (x, \ee)\leq v \leq V_{p(i+1)} (x, \ee) \right\} .
\end{equation}
Set 
$$ w_0 (k,x,\ee)= \max_{0\leq i \leq k} q_i \quad , \quad 
w_1 (k, x, \ee )= \max_{0\leq i \leq k} \left| \,  \left| 
V_{p(i+1)}(x, \ee) \right| - \left| V_{p(i)} (x, \ee) \right| \, \right| , $$
and 
$$ \Delta (k, x, \ee) = \max_{v \in \trr_{\Tileex} (\tileex)} \dist 
\left( v \;, \; \{ \varnothing ,  V_1 (x, \ee) , \ldots  , V_k (x, \ee) \} \right)  . $$
(\ref{rewriteq}) easily implies 
\begin{equation}
\label{boundmaxq}
w_0 (k,x,\ee) \leq  w_1(k,x,\ee)  + \Delta (k, x, \ee) .
\end{equation}
Since $\trr $ preserves height, we get for any $i\geq 1$,
$$ \left| V_i (x, \ee) \right|= \left| U_i (x, \ee) \right| = H_{\lfloor \u_i \# \tileex \rfloor } (\tileex  ) .$$
Then by Lemma \ref{convtileex} we get the following convergence in distribution 
$$ \ee \, w_1 (k,x, \ee)  \xrightarrow[\ee\rightarrow 0]{\quad  \quad} \max_{0\leq i \leq k} \left| D^{(x)}_{\u_{(i+1)}\zeta_x} -D^{(x)}_{\u_{(i)}\zeta_x }  
\right| \, ,$$
where $0=\u_{(0)} \leq \u_{(1)} \leq  \ldots \leq \u_{(k)} \leq \u_{(k+1)}=1 $ denotes 
the increasing re-ordering of $\{ 0, 1, \u_1 , \ldots  , \u_k \}$. Thus, 
\begin{equation}
\label{limitw1}
\forall \eta >0 \; , \quad \lim_{k\rightarrow \infty } \limsup_{\ee
  \rightarrow 0} \P \left( \ee \, w_1 (k, x, \ee ) > \eta \right) \, = \, 0.
\end{equation}
We next want to prove 
\begin{equation}
\label{limitdelta}
\forall \eta >0 \; , \quad \lim_{k\rightarrow \infty } \limsup_{\ee \rightarrow 0} \P \left( \ee \Delta (k, x, \ee ) > \eta \right) \, = \, 0.
\end{equation}
To that end, observe that for any $u, u'\in \tileex $ ,  $\dist ( \trr_{\Tileex}
(u) , 
\trr_{\Tileex} (u')) \leq \dist (u,u') $. Then if we set 
$$ \Delta' (k, x, \ee ) = \max_{u\in \tileex} \; \dist \left( u\; ; \;  \{ \varnothing , U_1 (x, \ee) , \ldots , U_k (x, \ee)  \} \right),  $$
we get 
\begin{equation}
\label{inegdelta}
\Delta (k, x, \ee) \leq \Delta' (k, x, \ee) 
\end{equation}
and we control $\Delta' (k, x, \ee) $ thanks to Lemma \ref{convtileex} (the
following 
argument is directly inspired 
from the proof of Theorem 20 \cite{Al2}): With
any $l\in \{ 0, \ldots , \# \tileex -1 \} $ we associate the index $i(l)\in \{ 0, \ldots , k+1 \} $ such that $\u_{i(l)}$ is the smallest element 
$y \in \{ 0,1, \u_1 , \ldots , \u_k \}$ such that $l\leq \lfloor y \# \tileex \rfloor $. Check that 
$$ \Delta' (k, x, \ee) \leq \max_{0\leq l < \# \tileex}  \left(  H_l (\tileex) +H_{\lfloor \u_{i(l)} \# \tileex \rfloor} (\tileex )
-2 \inf_{l\leq j \leq\lfloor \u_{i(l)} \# \tileex \rfloor} H_j (\tileex )  \right) .$$
Lemma \ref{convtileex} implies that the right member of the previous inequality converges in distribution to 
\begin{equation}
\label{quantity}
 \sup_{0\leq s \leq \zeta_x } \left(  D^{(x)}_s + D^{(x)}_{\u_{i(s)} \zeta_x}-2\inf_{s\leq  r \leq \u_{i(s)} \zeta_x} D^{(x)}_r \right) ,
\end{equation}
where we denote by $\u_{i(s)}\zeta_x $ the smallest element $y \in \{ \zeta_x,
\u_1 \zeta_x , \ldots , \u_k \zeta_x \}$ such that $s 
\leq y $ (recall that $\zeta_x$ stands for the lifetime of the process
$D^{(x)}$ as defined before Proposition \ref{theolimfini}). 
We easily check that (\ref{quantity}) converges to $0$ in
probability when $k$ goes to infinity since 
$$ \sup_{0 \leq s\leq \zeta_x} \left( \zeta_x \u_{i(s)} -s \right) \leq \max_{0 \leq i \leq k } \u_{(i+1)} - \u_{(i)} 
\xrightarrow[k\rightarrow \infty ]{\quad \quad} 0 $$
in probability. Thus, it implies (\ref{limitdelta}) by (\ref{inegdelta}) . Finally, as a consequence of (\ref{boundmaxq}),
(\ref{limitw1}) and (\ref{limitdelta}) we get 
\begin{equation}
\label{limitw0}
\forall \eta >0 \; , \quad \lim_{k\rightarrow \infty } \limsup_{\ee
  \rightarrow 0} \P \left( \ee \, w_0 (k, x, \ee ) > \eta \right) \, = \, 0.
\end{equation}
Then, check that on the event 
$$ E (k, x, \ee ,\delta) = \left\{  \, \min_{0 \leq i \leq k} \ee^2 \left( \Vbar_{p(i+1)} (x, \ee) -\Vbar_{p(i)} (x, \ee)  
\right) \; > \, \delta  \,  \right\} $$
the following inequality holds a.s.:
\begin{equation}
\label{omegw0}
\omega \left( H(x, \ee) , \delta \right) \leq 3 \, w_0 (k, x, \ee).
\end{equation}
Use Lemma \ref{convproba} to get 
$$\min_{0 \leq i \leq k} \ee^2 \left(  \Vbar_{p(i)} (x, \ee)  -\Vbar_{p(i+1)} (x, \ee) 
\right) \xrightarrow[\ee \rightarrow 0 ]{\quad  \quad} \frac{\zeta_x}{2\gamma}
\min_{0 \leq i \leq k} \left( \u_{(i+1)} - 
\u_{(i)} \right) $$
in distribution. Thus, 
$$\forall k \geq 1 \; , \quad \lim_{\delta \rightarrow 0} \liminf_{\ee
  \rightarrow 0} \P \left( E (k,x,\ee,\delta ) \, \right) = 1 $$
Easy arguments combined with (\ref{limitw0}) and (\ref{omegw0}) achieve the
proof of (T2) and at the same time the tightness for $H(x, \epsilon)$,
$\epsilon >0$. \cqfd

\vspace{3mm}

It remains to prove that $(2D^{(x)}_{\gamma s} \, ; s\geq 0)$ is the only
possible weak limit for the processes $H(x, \epsilon)$,
$\epsilon >0$.
Tightness for the $H(x, \epsilon)$'s, $\epsilon >0$ and Lemma 3.6 imply 
that the joint distributions of 
$(H(x, \ee ) \; , \; 2\ee^2 \# \tileex )$ , $\epsilon >0$ are tight. 
Assume that along a subsequence 
$\ee_p \rightarrow 0$ the following joint convergence 
$$ \left( H(x, \ee_p ) \; , \; 2\ee_p^2 \# \tileexp \right)
\xrightarrow[p \rightarrow \infty ]{\quad {\rm d} \quad} \left( H' , 
\zeta' \right) $$
holds for some continuous 
process $H'$ and some positive random variable $\zeta'$. 
Lemma \ref{convproba} implies 
\begin{multline*}
\left(  H(x, \ee_p ) \; ; 2\ee_p^2 \# \tileexp \, ; 2
\ee_p^2 \Vbar_1 (x, \ee_p)
  \, , \ldots , 2 \ee_p^2 \Vbar_k (x , \ee_p) \right) \\ 
\xrightarrow[\ee\rightarrow 0]{\quad {\rm d} \quad} \left( H' \, ;  
\zeta'  \, ; \u_1 \zeta'/ \gamma \, , \ldots , \, 
\u_k \zeta' / \gamma \right) 
\end{multline*}
where the $\u_i $'s are chosen independent of $(H' , \zeta')$. 
Since $\trr $ preserves height, we get for any 
$i \geq 1 $
$$ H_{ \Vbar_i (x, \ee_p )} (\txe ) = \left| V_i (x, \ee_p ) \right| = \left| U_i (x, \ee_p ) \right| = H_{\lfloor 
\u_i \# \tileexp \rfloor} (\tileexp ) .$$
Then Lemmas \ref{convtileex} and \ref{convproba} imply for any $k\geq 1$ 
\begin{multline*}
\left( H_{2\ee_p^2 \Vbar_1 (x, \ee_p )} (x, \ee_p ) \, , \, \ldots , \, 
H_{2\ee_p^2 \Vbar_k (x, \ee_p )} (x, \ee_p ) 
\, ; \, 2 \ee_p^2 \# \tileexp \, ; \,  
2\ee_p^2 \Vbar_1 (x, \ee_p )\, , \ldots , \, 2\ee_p^2 \Vbar_k (x, \ee_p )
\right) \\
\xrightarrow[\ee\rightarrow 0]{\quad  \quad} 
\left( 2D^{(x)}_{\u_1 \zeta_x} \, , \ldots , \,  2D^{(x)}_{\u_k \zeta_x} \, ; \,
  \zeta_x \, ; \, 
\u_1 \zeta_x / \gamma \, , \ldots , \,\u_k \zeta_x / \gamma \right) 
\end{multline*}
in distribution. Consequently, 
\begin{multline*}
\left( H'_{\u_1 \zeta' /\gamma } \, , \ldots , \, H'_{\u_k \zeta'/\gamma } \,
  ; \, 
\zeta' \, ; \, \u_1\zeta' /\gamma \, , \ldots , \, \u_k \zeta' /\gamma \right)
\\ \laweq   
\left( 2D^{(x)}_{\u_1 \zeta_x} \, , \ldots , \, 2D^{(x)}_{\u_k \zeta_x} \, ; \, \zeta_x \, ; \, \u_1\zeta_x /\gamma \, , \ldots , \, 
\u_k \zeta_x /\gamma \right) .
\end{multline*}
It implies $(H_s' \, ; s\geq 0) \laweq (2D^{(x)}_{\gamma s} \, ; s\geq
0)$, which achieves the proof of 
(\ref{mainconv}). \cqfd

\subsection{Proof of Lemma \ref{convproba} }
We introduce the notation 
$$\Ubar_i (x,\ee)= \sum_{\substack{v\in \U  \\ v\leq V_i(x, \ee)}} Z_v (\Tileex )   $$
and we first prove the following convergence in probability:
\begin{equation}
\label{convubar}
\ee^2 \left( \Ubar_i (x, \ee) - \u_i \# \tileex \right)
\xrightarrow[\ee\rightarrow 0]{\quad  \quad} 0.
\end{equation}

\noindent 
{\bf Proof:}
Deduce from (\ref{combifond2}): 
$$\{ u\in \tileex :  \trr_{\Tileex } (u) < V_i (x, \ee) \}   \subset \{
u\in \tileex :  u \leq  U_i (x, \ee) \} 
\subset \{  u\in \tileex :  \trr_{\Tileex } (u) \leq  V_i (x, \ee)\} $$
which implies 
\begin{equation}
\label{coinc}
0 \leq \Ubar_i (x, \ee) - \lfloor \u_i \# \tileex \rfloor \leq Z_{V_i (x, \ee )} (\Tileex ). 
\end{equation}
Then observe that for any $v \in \U $, 
$$ \P \left( V_i (x, \ee ) =v \left| \Tileex  \right. \right) = \frac{Z_v (\Tileex ) }{\# \tileex} .$$  
Thus, (\ref{coinc}) and Cauchy-Schwarz inequality imply 
\begin{eqnarray*}
\E \left[ \left|  \Ubar_i (x, \ee) - \u_i \# \tileex \right| \right] &\leq & 1+ \E \left[ \frac{1}{\# \tileex} \, \sum_{v\in \U}  Z_v(\Tileex)^2 
\right] \\
&\leq & 1+ \E \left[ \frac{1}{(\# \tileex)^2} \, \sum_{v\in \U}
  Z_v(\Tileex)^2 \right]^{1/2}\E \left[\sum_{v\in \U}
  Z_v(\Tileex)^2 \right]^{1/2}. 
\end{eqnarray*}
Since $\# \tileex = \sum_{v \in \U} Z_v (\Tileex ) $, we get $\sum_{v\in \U}
  Z_v(\Tileex)^2  \leq  (\# \tileex)^2$ and  
\begin{equation}
\label{coincl1}
\E \left[ \left|  \Ubar_i (x, \ee) - \u_i \# \tileex \right| \right] \leq 1+ \E \left[\sum_{v\in \U}  Z_v(\Tileex)^2 
\right]^{1/2} .
\end{equation}
Remark \ref{shuffleGWI}
and (\ref{combifond1}) imply that $\tilee$ is a GWI-tree with immigration distribution $\nu =\mu $, so that 
$$g(x)= f(x)= \frac{u}{1-dx} \quad {\rm and } \quad g^{(j)}(1)= j! \left( \du \right)^j \; , \; j\geq 1 .$$ 
Then, by Proposition \ref{factmomentZGWI} $(iii)$
$$ \E \left[\sum_{v\in \U}  Z_v(\Tileex)^2 
\right] \leq \frac{ K_{\aa }\, \xxe}{1-d/u} =  \frac{K_{\aa} x }{\ee^2} \, (1+ o(1)) .$$
Thus, 
$$ \E \left[ \ee^2 \left|  \Ubar_i (x, \ee) - \u_i \# \tileex \right| \right]
\leq K_{\aa , x } \, \ee (1 + o(1)) $$
and (\ref{convubar}) follows. \cqfd

\vspace{4mm}

Then, Lemma \ref{convproba} is a consequence of
the convergence in probability: 
\begin{equation}
\label{convvbar}
\ee^2 \left( \Vbar_i (x, \ee) - \frac{1}{\gamma} \, \Ubar_i (x, \ee) \right) \xrightarrow[\ee\rightarrow 0]{\quad  \quad} 0 .
\end{equation}

\noindent 
{\bf Proof of} (\ref{convvbar}): We need several preliminary
estimates (Lemmas \ref{betaGW} and \ref{spinall2}) whose proofs rely on Propositions \ref{ZestimateGW} and
\ref{factmomentZGWI}. We first consider a 
random marked GW-forest with $l$ elements 
$\f_{\ee }=(\varphi_{\ee} \, ; (\mu_u , u\in \varphi_{\ee}))$ as defined at
Proposition \ref{ZestimateGW}: recall that  $\varphi_{\ee} =(\tau_1 , \ldots
\tau_l)$ is a forest of $l$ i.i.d  GW($\mu $)-trees and that the marks $(\mu_u , u \in
\varphi)$ are i.i.d. conditional on $\varphi_{\ee}$, their conditional
distribution being given by $\a$. Set 
$\t_{1,\ee} =(\tau_1 \, ; (\mu_u , u\in \tau_1))$ and define 
$$ 1/ \gamma_{\ee} = \frac{ \E \left[  \sum_{v\in \U }  \un_{\{ Z_v (\t_{1, \ee}) >0
      \} }\right] }{\E \left[  \sum_{v\in \U }  Z_v (\t_{1, \ee} ) \right]} .$$
We also set 
$$ \beta \left( \f_{\ee } \right) =\sum_{v\in \U }  Z_v
(\f_{\ee })- \gamma_{\ee} \, \un_{\{ Z_v (\f_{\ee }) >0
      \} } . $$
\begin{lemma}
\label{betaGW} 
First $(i)$ $\lim_{\ee \rightarrow 0} \gamma_{\ee}= \gamma \;  $ 
and for any $l\geq 1$, 
$$ (ii) \, \quad \quad 0\leq \E \left[ \beta  \left( \f_{\ee } \right) \right] \leq
K_{\aa } \, l(l-1) , $$
$$ (iii) \,  \quad \quad  \E \left[ \beta  \left( \f_{\ee } \right)^2 \right] \leq
K_{\aa } \, \frac{l^4}{1-d/u} \quad {\rm and }\; {\rm thus} \quad 
\E \left[ \left| \beta  \left( \f_{\ee } \right) \right| \right] \leq
K_{\aa } l^2 (1-d/u)^{-1/2} .$$
\end{lemma}
\noindent {\bf Proof :} Let us prove $(i)$: First observe that
$$ 1/ \gamma_{\ee} = 
\frac{  \sum_{v\in \U }  1-f_v (0)  }{\sum_{v\in \U }f_v'(1) } .$$
Then, Proposition \ref{ZestimateGW} implies that for any $v=m_1 \ldots m_n \in
\U$
$$ \sum_{v\in \U } f_v'(1)= 
\sum_{v\in \U } a_v (d/u)^{|v|} 
=\frac{1}{1-d/u} $$
and also 
$$  1-f_v (0) = a_v \left( \du \right)^{|v|} \left[  \left( \du \right)^{n} a_{m_n} \ldots a_{m_1}
+ \left( \du \right)^{n-1} a_{m_n} \ldots a_{m_2} + \ldots +1 \right]^{-1} ,$$
Thus, we get  
$$ 1/ \gamma_{\ee} 
= \E  \left[  \left( 1+ X_1 \du + \ldots + X_1 X_2 \ldots X_G
    \left( \du \right)^{G}  \right)^{-1} \right] . $$
where we recall that the sequence of random variables 
$(X_n ; n\geq 0)$ is distributed as
specified after formula (\ref{coeff}), and 
where $G$ stands for an independent random variable whose distribution is 
given by
$\P (G=n)=(1-d/u)(d/u)^n $, $n\geq 0$. Since $\lim_{\ee \rightarrow 0} d/u
=1$, an elementary argument implies
$$ \lim_{\ee \rightarrow 0} 1/\gamma_{\ee} = \E \left[ \left( 1+ X_1 + X_1 X_2 + X_1 X_2 X_3 + \ldots
  \right)^{-1} \right] =1/ \gamma .$$

Let us prove $(ii)$: Deduce from Proposition \ref{ZestimateGW} that 
\begin{equation}
\label{precedente}
\E \left[ \beta  \left( \f_{\ee } \right) \right] = \sum_{v\in \U } lf_v'
(1) - \gamma_{\ee} \, \left( 1-f_v (0)^l \right).
\end{equation}
The definition of $\gamma_{\ee}$ implies 
$$ l\sum_{v\in \U } f_v'
(1) - \gamma_{\ee} ( 1-f_v (0) ) =\sum_{v\in \U } lf_v'
(1) - l\gamma_{\ee} ( 1-f_v (0) ) =0 .$$
We then subtract this expression to (\ref{precedente}) and we get
$$\E \left[ \beta  \left( \f_{\ee } \right) \right] = \gamma_{\ee}
\sum_{v\in \U } f_v (0)^l -1 + l \left( 1-f_v (0) \right) .$$ 
Then, use the elementary inequality $(1-x)^l-1+l x \leq l(l-1)x^2/2$,  $x\in
[0,1]$ to get 
\begin{equation}
\label{lunbeta1} 
\E \left[ \beta  \left( \f_{\ee } \right) \right] \leq  \frac{\gamma_{\ee}l(l-1)}{2}\sum_{v\in \U } \left(
  1-f_v (0) \right)^2 
\end{equation}
Deduce from the explicit computation of $1-f_v
(0)$ recalled above that 
$$ \left( 1-f_v (0) \right)^2 \leq 
\left( \du\right)^{2 |v| } a_v^2 \leq
a_+^{ |v|}  a_v .$$
Thus, 
$$ \sum_{v\in \U } \left( 1-f_v (0) \right)^2  \leq \sum_{n\geq 0} 
a_+^{ n} \sum_{ m_1, \ldots , m_n \in \N^* }
a_{m_1} \ldots a_{m_n} \leq (1-a_+)^{-1} $$
and $(ii)$ follows from $(i)$. 

\vspace{3mm}

  It remains to prove $(iii)$: For convenience of notation, we simply write 
$\beta $ and $Z_v$ instead of $\beta  \left( \f_{\ee } \right) $ and $Z_v
\left( \f_{\ee } \right)$. Check that 
\begin{equation}
\label{betacarre}
\E \left[ \beta^2 \right] = \E \left[ E_1 \right]+ \E \left[ E_2 \right], 
\end{equation}
where 
$$  E_1 = \sum_{\substack{v,v'\in \U \\  v\wedge v'\notin\{ v,v'\} }} \left(
 Z_v - \gamma_{\ee} \, \un_{\{ Z_v >0\} } \right) \left(
 Z_{v'} - \gamma_{\ee} \, \un_{\{ Z_{v'} >0\} } \right) $$
and 
$$  E_2 = \sum_{\substack{v,v'\in \U \\  v\wedge v'\in\{ v,v'\} }} \left(
 Z_v - \gamma_{\ee} \, \un_{\{ Z_v >0\} } \right) \left(
 Z_{v'} - \gamma_{\ee} \, \un_{\{ Z_{v'} >0\} } \right) $$
(note that in the two sums all but a finite number of terms vanish). Define for any $w\in \U$
$$ \beta_w = \sum_{v\in \U }  Z_{wv} - \gamma_{\ee} \, \un_{\{ Z_{wv} >0\} }.$$
$E_1$ can be rewritten as follows
$$ E_1 = \sum_{w\in \U } \sum_{i \neq j\in \N^*} \beta_{wi}\beta_{wj}.$$
Deduce from Proposition \ref{ZestimateGW} that conditional on 
$(Z_{wi}, Z_{wj})$ (with $i\neq j$) the random variables 
$\beta_{w_i}$ and $\beta_{w_j}$ are independent and
distributed as $\beta$ with resp. $l=Z_{wi}$ and $l=Z_{wi}$. Use $(ii)$ to get 
\begin{eqnarray*}
\E \left[ \beta_{w_i}\beta_{w_j} \mid (Z_{wi}, Z_{wj}) \right] &=& \E \left[
  \beta_{w_i} \mid Z_{wi}\right]
\E \left[ \beta_{w_j} \mid Z_{wj}\right] \\
&\leq & K_{\aa} Z_{wi}(Z_{wi}-1) Z_{wj}(Z_{wj}-1).
\end{eqnarray*}
By Proposition \ref{ZestimateGW} again, we get  
$$ \E \left[ x^{ Z_{wi} } y^{ Z_{wj} } \mid Z_w \right]= f \left( 1-a_i-a_j
  +a_i x+a_j y\right)^{Z_w} .$$
Recall that
\begin{equation}
\label{derivf}
\frac{\der^k f^n}{\der x^k} (x)= \left( \du \right)^k \frac{(n+k-1)!}{(n-1) !}
f(x)^{n+k}.
\end{equation}
Then, 
\begin{equation*}
 \E \left[ Z_{wi}(Z_{wi}-1) Z_{wj}(Z_{wj}-1) \mid Z_w \right] = a_i^2  a_j^2
 \left( \du \right)^4 \left( Z_w +3\right)_4 
\leq  12 Z_w^4.
\end{equation*}
Thus, by Proposition \ref{ZestimateGW}
\begin{equation}
\label{eun2}
\E \left[ E_1 \right] \leq  K_{\aa } \, \E \left[ \sum_{w\in \U }Z_w^4
\right] \leq  K_{\aa } \, \frac{l^4}{1-d/u} .
\end{equation}

We get a similar upper-bound for $\E \left[ E_2 \right]$ by first noting that 
$$ E_2 \leq 2 \sum_{w\in \U } \left( Z_w - \gamma_{\ee} \, \un_{\{ Z_w
    >0\} } \right) \beta_{w} .$$
Apply Proposition \ref{ZestimateGW} $(i)$ and Lemma \ref{betaGW} $(ii)$ to get 
$$ \E \left[ \left( Z_w - \gamma_{\ee} \, \un_{\{ Z_w >0\} } \right)
  \beta_{w} \mid Z_w \right] \leq K_{\aa} Z_w^3 .$$
By Proposition \ref{ZestimateGW} $(iii)$ again 
\begin{equation}
\label{edeux}
\E \left[ E_2 \right] \leq  K_{\aa } \, \E \left[ \sum_{w\in \U }Z_w^3
\right] \leq  K_{\aa }\frac{l^3}{1-d/u} .
\end{equation}
Then $(iii)$ follows from (\ref{betacarre}), (\ref{eun2}) and
(\ref{edeux}). \cqfd 

\vspace{4mm}

We need similar estimates for a marked GWI($\mu , r$)-forest $\f_{0, \ee}$ 
whose distribution is the same as in Proposition \ref{factmomentZGWI}: recall
that $r$ is some fixed repartition
probability measure on $\{ (k,l)\in \N^* \times \N^* \, :\, l\leq k \}$. We denote by $\nu$ the corresponding 
immigration probability measure given by $\nu (k-1)=\sum_{1\leq l\leq k}
r(k,l)$ , $k\geq 1$ and we set  $g(r)= \sum_{k\geq 0} \nu (k)r^k $. We define
$\f_{0, \ee}$ as $(\varphi_{0, \ee} ; (\mu_u , u \in \varphi_{0, \ee} ))$ 
where $\varphi_{0, \ee} =
(\tau_0 , \tau_1, \ldots , \tau_l)$, the $\tau_i$'s are mutually independent,   
$\tau_1 , \ldots , \tau_l$ are i.i.d. GW($\mu $)-trees,  $\tau_0$ is a GWI($\mu
, r$)-tree and conditional on $\varphi_{0, \ee}$ 
the marks $\mu_u$ are i.i.d. 
random variables distributed in accordance with $\a$. Recall notations 
$$ u_n^*= u^*_n(\varphi_{0, \ee}) \; ,\;   v^*_n = \trr_{\f_{0, \ee}} (u^*_n)
\; , \quad \spin=\{ v^*_n i \, , \; i\in \N^* \setminus \{ \mu_{u^*_n}\}\, , n \geq 0 \} $$ 
and recall that $\s $ is the $\sigma$-field generated 
by the random variables  $(\mu_{u^*_n} \, ; n\geq 0)$ and $(Z_w (\f_{0, \ee})
\; , w\in \spin ) $. For any $n\geq 1$ we also set 
$$ \spin (n)= \{ w\in \spin \, : |w| \leq
n\} \cup \{  v^*_n \} .$$
We set 
$$ \beta_w \left( [\f_{0, \ee}]_{ u_n^*} \right)  = 
\sum_{v\in \U } Z_{wv} (  [\f_{0, \ee}]_{ u_n^*} ) - \gamma_{\ee}
\, \un_{\{ Z_{wv} ( [\f_{0, \ee}]_{ u_n^*} )>0\} } .$$
\begin{lemma}
\label{spinall2}
For any $n\geq 1$,
$$ \E \left[ \sup_{A\subset \spin (n)} \left| 
\sum_{w\in A}\beta_w ( [\f_{0, \ee}]_{ u_n^*} )
  \right| \right] \leq 
K_{\aa} \, n (l+1)^2 (1-d/u)^{-1/2} 
\, \max \left( 1, g'(1)^2, g''(1)^2 \right). $$
\end{lemma}
\noindent {\bf Proof :} To simplify notation we write $\beta_w $ and $Z_w$
instead of $\beta_w ([\f_{0, \ee}]_{ u_n^*} )$ and  $Z_{w} (
[\f_{0, \ee}]_{ u_n^*})$. We also
denote by $\E^{\s}$ the $\s$-conditional expectation. Let $A\subset
\spin(n)$. From Proposition \ref{factmomentZGWI} $(i)$ 
we deduce that conditional
on $\s $ the $(\beta_w\, ; w\in \spin (n))$ are independent random variables 
and that for each $w\in \spin (n)$, conditional on $Z_w=l$,     
$\beta_w $ is distributed as the random variable 
$\beta (\f_{\ee})$ defined at Lemma \ref{betaGW}. 
Apply Lemma \ref{betaGW} to get
\begin{eqnarray*}
\E^{\s} \left[\left| \sum_{w\in A}\beta_w 
  \right| \right] &\leq & \sum_{w\in A} \E^{\s} \left[ \left| \beta_w \right|
  \right] \\
& \leq & K_{\aa} (1-d/u)^{-1/2} \sum_{w\in \spin (n)} Z_w^2 .
\end{eqnarray*}
Next, use Proposition \ref{factmomentZGWI} $(ii)$ to get 
\begin{eqnarray*}
\E \left[ \sum_{w\in \spin (n)}  Z_w^2 \right] & \leq & 
\sum_{k=0}^{n-1} 
\sum_{i\in \N^*} \E \left[ Z_{v^*_k i}^2 \right] \\
&\leq & K_{\aa} \,  n\, (l+1)^2 \, \max \left( 1, g'(1)^2,  g''(1)^2\right),
\end{eqnarray*}
which achieves the proof of the lemma. \cqfd

\vspace{4mm}

We now comme back to the proof of (\ref{convvbar}) and we apply the previous
results to the marked sin-tree $\Tilee = (\tilee ; (\muti_u, u\in
\tilee))$. For convenience of notation, we fix $i$ and we set 
$$ U=U_i (x,\ee)\; , \;  V=V_i (x,\ee) \; ,\; \Ubar=\Ubar_i (x,\ee) \; , \; 
\Vbar=\Vbar_i (x,\ee) .$$
We keep the notations $u_n^*= u^*_n(\tilee)$, $v^*_n = \trr_{\Tilee}
(u^*_n)$, $\spin$, $\spin (n) $ and  $\s $. Recall that $\Tileex= [\Tilee]_{u^*_{\xxe}}$ and that for any
$v\in \U$ that is not a descendant of $v^*_{\xxe} $ 
\begin{equation}
\label{cutuncut}
Z_w \left( \Tileex \right)=  Z_w \left( \Tilee \right).
\end{equation}
For convenience of notation, we set for any $w\in \U$ 
$$ Z_w = Z_w \left( \Tileex \right) \quad {\rm and } \quad \beta_w =
\sum_{v\in \U} Z_{wv}-\gamma_{\ee} \un{ \{ Z_{wv} >0\} } .$$
Since $\lim_{\ee \rightarrow 0}\gamma_{\ee} =\gamma $, we only have to show 
\begin{equation}
\label{convbabar}
\ee^2 \left(  \Ubar - \gamma_{\ee} \Vbar \right) = \ee^2 \sum_{v\leq V } Z_v - \gamma_{\ee} \un_{\{ Z_v >0 \}}  
\xrightarrow[\ee\rightarrow 0]{\quad  \quad} 0 
\end{equation}
in probability. To that end, we first introduce the random word 
$$ W= \max \{w \in \spin(\xxe) \, : \; w \leq V \} ,$$ 
where the maximum is taken with respect
to the lexicographical order on $\U$. There are two cases:

\noindent
$\bullet $ If $V\notin \llb \varnothing , v^*_{\xxe -1}\rrb $, then we can find
$V' \in \U $ such that $V=WV'$ and we set in that case $A=\{ w \in
\spin(\xxe) \, : \; w <W \}$. 

\noindent
$\bullet $ If $V \in \llb \varnothing , v^*_{\xxe -1}\rrb $, then we set 
$A=\{ w \in
\spin(\xxe) \, : \; w \leq W \}$.

\vspace{3mm}

\noindent
Then, check that $ \Ubar - \gamma_{\ee} \Vbar = e_1 (\ee) +e_2 (\ee) +e_3 (\ee) $ with
$$ e_1 (\ee)= \sum_{\substack{v\in \llb \varnothing , v^*_{\xxe -1}\rrb \\ v\leq V}} Z_{v}- \gamma_{\ee} \un_{\{ Z_{v} >0  \}} ,$$

$$ e_2 (\ee)= \sum_{w\in A} \beta_w \quad {\rm and }\quad 
e_3 (\ee) = \un_{\{V\notin \llb \varnothing , v^*_{\xxe -1}\rrb  \}} \sum_{v \leq V' } Z_{Wv}- \gamma_{\ee} \un_{\{ Z_{Wv} > 0 \}} .$$
The limit (\ref{convbabar}) is then implied by the following convergences
\begin{equation}
\label{convere1}
\ee^2 \E \left[ |e_1 (\ee ) | \right] \xrightarrow[\ee\rightarrow 0]{\quad  \quad} 0 \; ,
\end{equation}

\begin{equation}
\label{convere2}
\ee^2 \E \left[ \left| e_2 (\ee ) \right| \right] \xrightarrow[\ee\rightarrow 0]{\quad  \quad} 0 \; ,
\end{equation}

\begin{equation}
\label{convere3}
\ee^4 \E \left[ e_3 (\ee )^2 \right] \xrightarrow[\ee\rightarrow 0]{\quad  \quad} 0\; .
\end{equation}

\noindent {\bf Proof of (\ref{convere1}) :} Use Proposition \ref{factmomentZGWI} $(ii)$ 
with $p=1$, $l=0$, $n=\xxe -1$ and $g(x)=f(x)=u/(1-dx)$ to get 
\begin{eqnarray*}
\E \left[ |e_1 (\ee ) | \right] & \leq & \sum_{i=0}^{\xxe -1} \E \left[
  Z_{v^*_i}\right] + (\xxe -1)\gamma_{\ee}  \\
&\leq &  K_{\aa} d (\xxe -1) /u  + (\xxe -1) \gamma_{\ee} \leq K_{\aa , x} \, \ee^{-1}
\end{eqnarray*}
which obviously implies (\ref{convere1}). \cqfd

\vspace{4mm}

\noindent {\bf Proof of (\ref{convere2}) :} We use Lemma \ref{spinall2} with
$n=\xxe $, $l=0$ and $g(x)=f(x)=u/(1-dx)$ 
and thus $g^{(j)}(1)= j!(d/u)^j $, to get 
$$\E \left[ \left| e_2 (\ee ) \right|   \right] \leq K_{\aa} \xxe
(1-d/u)^{-1/2} \leq K_{\aa ,x} \, \ee^{-3/2} $$
which implies (\ref{convere2}). \cqfd 

\vspace{4mm}

\noindent {\bf Proof of (\ref{convere3}) :} It requires more complicated 
arguments. Let $w_0 \in \U $ and let $l$ be a positive integer. We define $E(w_0
,l)$ as the event 
$\{ W=w_0 \; ; \; Z_{w_0} =l \}$. We first get an upper-bound for 
$$ \xi (w_0 , l)= \E \left[ e_3 (\ee)^2 \mid E (w_0, l) \right]. $$
Let $\f=(\varphi\, ; (\mu_u , u\in
\varphi))$ be a marked GW-forest with $l$ elements as defined 
at Proposition \ref{ZestimateGW}. Pick uniformly at random a vertex $\u (\f)$ in
$\varphi $ and define $\v (\f)\in \U$ by $\v (\f)=\trr_\f (\u (\f))$. As a consequence
of Propositions \ref{ZestimateGW} $(i)$ and \ref{factmomentZGWI} $(i)$, we get the following identity
\begin{equation}
\label{identlaw}
\left( Z_{w_0 v} \; ,\; v\in \U \; ; \; V' \right) \quad {\rm under} \quad
\P \left( \; \cdot \mid E (w_0, l)\right) 
\; \laweq \; \left( Z_v (\f) \; , v\in \U \; ; \; \v (\f) \right).
\end{equation}
Let $G$ be the function on $\F $ defined by 
$$ G\left( [\f]_{\u (\f)} \right) = \sum_{v \leq \v (\f) } Z_v (\f) - \gamma_{\ee} \un_{\{ Z_v (\f) >0 \}}. $$ 
Then, (\ref{identlaw}) implies 
\begin{eqnarray}
\label{ineqgrossun}
\xi (w_0 , l) &=& \E \left[  G\left( [\f]_{\u (\f)} \right)^2 \right]= \E \left[ \frac{1}{\# \varphi} 
\sum_{u\in \varphi} G\left( [\f]_{u} \right)^2 \right] \\
\label{ineqgrossdeux} 
&\leq & (1+\gamma_{\ee}) \E \left[ \sum_{u\in \varphi} \left| G\left( [\f]_{u}
    \right) \right| \right] , 
\end{eqnarray}
since for any $u\in \varphi $, 
$$ \frac{1}{\# \varphi} \left| G\left( [\f]_{u}  \right) \right| \leq
\frac{1+\gamma_{\ee }}{\# \varphi} \sum_{v \leq \trr_\f (u) } Z_v (\f) \leq
1+\gamma_{\ee } .$$

We now estimate the right member of (\ref{ineqgrossdeux}) thanks to (\ref{sizebias}): 
Recall the notation $\varphi_{\flat}$ for a size-biased forest with $l$ elements, i.e. a 
GWI($\mu,r$)-forest with $l$ elements where $r$ is given by  
$r(k,j)= ud^{k} / \bar{\mu}$ , $ 1\leq j \leq k $ with $\bar{\mu}= \sum_{k\geq 0} k\mu(k)=d/u$.
Thus the corresponding immigration distribution is $\nu (k)= (k+1) u^2 d^{k}$
,  $k\geq 0$ and its generating function is 
$g(r)=u^2/(1-dr)^2 $. Let us define the random marked $GWI$-forest
$\f_{\flat} $ as $ (\varphi_{\flat} \, ; (\mu_u^{\flat} \in \varphi_{\flat}))$ where conditional on $\varphi_{\flat}$ the $\mu^{\flat}_u$'s are i.i.d. with 
distribution ${\bf a}$. 
Deduce from (\ref{sizebias}) that 
\begin{equation}
\label{applysizebias}
\E \left[ \sum_{u\in \varphi}  \left| G\left( [\f]_{u} \right) \right| \right] = \sum_{n\geq
  0} l \left( \du \right)^n \E \left[ \left| G \left([\f_{\flat}]_{u^*_n
  (\f_{\flat})} \right)  \right|  \right] .
\end{equation} 

Set as usual $v_n^* (\f_{\flat}) = {\bf Tr}_{\f_{\flat}} (u_n^* (\f_{\flat}))$
and observe for any $n\geq 0 $  
\begin{eqnarray*}
G \left([\f_{\flat}]_{u^*_n (\f_{\flat})} \right) &=& \sum_{v < v^*_n (\f_{\flat} )} Z_v (\f_{\flat}) -\gamma_{\ee} \un_{\{ Z_v (\f_{\flat}) >0 \} } \\
&=& \sum_{w\in A_{\flat}} \beta_w (\f_{\flat} ) + \sum_{i=0}^{n-1}  Z_{v^*_i} (\f_{\flat}) -\gamma_{\ee} 
\end{eqnarray*}
where we have set 
$$ \beta_w (\f_{\flat}) = \sum_{v \in \U} Z_{wv} (\f_{\flat} ) -\gamma_{\ee} \un_{\{ Z_{wv} (\f_{\flat}) >0 \}} $$
and $A_{\flat}= \{ w\in \spin_{\flat}  \, : \; w< v^*_n (\f_{\flat}) \}$ with 
$$\spin_{\flat}= \{ v^*_{k-1} (\f_{\flat}) i \; ;\;  i\in \N^*
\setminus \{ 
\mu_{u^*_{k}(\f_{\flat})} \} \, , \, k\geq 1 \} .$$
Then, 
$$ \E \left[ \left| G \left([\f_{\flat}]_{u^*_n
  (\f_{\flat})} \right)  \right|  \right]  \leq \E \left[ \left| \sum_{w\in
A_{\flat}} \beta_w (\f_{\flat} ) \right| \right] + \E \left[ \left| 
\sum_{i=0}^{n-1}  Z_{v^*_i} (\f_{\flat}) -\gamma_{\ee} 
\right| \right] . $$
Use Lemma \ref{spinall2} with 
$g(x)= u^2/(1-dx)^2 $ to get
$$ \E \left[ \left| \sum_{w\in
A_{\flat}} \beta_w (\f_{\flat} ) \right| \right]  \leq  K_{\aa} \, n l^2
(1-d/u)^{-1/2} $$
and use Proposition \ref{factmomentZGWI} $(ii)$ with $p=1$ and
$g(x)=u^2/(1-dx)^2$ to get 
$$ \E \left[ \left| 
\sum_{i=0}^{n-1}  Z_{v^*_i} (\f_{\flat}) -\gamma_{\ee}
\right| \right] \leq  K_{\aa} \, n l .$$
These inequalities imply 
\begin{equation}
\label{boundalpha}
\xi (w_0 , l) \leq K_{\aa} \, l^3 (1-d/u)^{-1/2} \sum_{n\geq 0} n\left( \du
\right)^n \leq K_{\aa} \, \frac{l^3}{(1-d/u)^{5/2}}. 
\end{equation}
We now comme back to the proof of (\ref{convere3}): by (\ref{boundalpha}), we
get 
\begin{eqnarray*}
\E \left[ e_3 (\ee)^2 \right] &=& \sum_{\substack{ w_0 \in \U ,\\ l \geq 1}} \xi (w_0 , l) 
\P(W=w_0 \; ; \; Z_{W}=l) \\
&\leq & \frac{K_{\aa}}{(1-d/u)^{5/2}} \E \left[ \sum_{w_0 \in \spin (\xxe) }
  Z_{w_0}^3 \un_{\{ W=w_0 \}}  \right] 
\end{eqnarray*}
Then, set for any $w_0\in \spin (\xxe) $, $\zeta_{w_0} = \sum_{v\in \U } Z_{w_0 v} $ and observe 
that $\P (W=w_0 \mid \s)= \zeta_{w_0}/ \#\tileex$. Thus the previous inequality implies
\begin{eqnarray*}
\E \left[ e_3 (\ee)^2 \right] &\leq &
\frac{K_{\aa}}{(1-d/u)^{5/2}} \E \left[ \sum_{w_0 \in \spin (\xxe) }  Z_{w_0}^3 \frac{\zeta_{w_0}}{\#\tileex} \right]  \\
 &\leq &
\frac{K_{\aa}}{(1-d/u)^{5/2}} \E \left[ \sum_{w_0 \in \spin (\xxe) }  Z_{w_0}^{6}  \right]^{1/2} 
\E \left[ \sum_{w_0 \in \spin (\xxe) }  \frac{\zeta_{w_0}^2}{(\# \tileex)^2}\right]^{1/2} 
\end{eqnarray*}
But $\sum_{w_0 \in \spin (\xxe) } \zeta_{w_0}^2 \leq (\# \tileex)^2$ since $1+
\xxe + \sum_{w_0 \in \spin (\xxe) } 
\zeta_{w_0} = \# \tileex $. Then, use Proposition \ref{factmomentZGWI} $(ii)$
with $p=6$, $l=0$ and $g(x)=f(x)$ to get  
$$ \E \left[ e_3 (\ee)^2 \right] \leq K_{\aa} \frac{ \xxe^{1/2}}{(1-d/u)^{5/2}}
\leq K_{\aa , x} \, \ee^{-3} , $$
which implies (\ref{convere3}). \cqfd

\bibliographystyle{plain}

\end{document}